\documentclass[12pt,twoside]{article}

\begin{document}

\date{}




\centerline{}

\centerline{}

\centerline {\Large{\bf A Note on Generalized Intuitionistic Fuzzy }}

\centerline{}

\centerline {\Large{\bf $\psi$ Normed Linear Space}}

\centerline{}


\centerline{\bf {Sumit Mohinta and T.K. Samanta}}

\centerline{}

\centerline{e-mail: sumit.mohinta@yahoo.com}

\centerline{Department of Mathematics, Uluberia College, India-711315.} %
\centerline{e-mail: mumpu$_{-}$tapas5@yahoo.co.in}

\begin{abstract}
\textbf{\emph{\ In respect of the definition of
intuitionistic fuzzy n-norm \cite{Vijayabalaji} ,
the definition of generalised intuitionistic fuzzy $\psi$ norm (\, in short GIF$\psi$N \,) is introduced over a linear space and there after a few results on generalized intuitionistic fuzzy $\psi$ normed linear space and finite dimensional generalized intuitionistic fuzzy $\psi$ normed linear space have been developed. Lastly, we have introduced the
definitions of generalised intuitionistic fuzzy $\psi$ continuity , sequentially
intuitionistic fuzzy $\psi$ continuity and it is proved that they are equivalent.}}
\end{abstract}

\centerline{}

\newtheorem{Theorem}{\quad Theorem}[section]

\newtheorem{Definition}[Theorem]{\quad Definition}

\newtheorem{Corollary}[Theorem]{\quad Corollary}

\newtheorem{Lemma}[Theorem]{\quad Lemma}

\newtheorem{Note}[Theorem]{\quad Note}

\newtheorem{Remark}[Theorem]{\quad Remark}

\newtheorem{Result}[Theorem]{\quad Result}

\newtheorem{Example}[Theorem]{\quad Example}




\section{Introduction}

In 1965, Zadeh\cite{zadeh} first introduced the concept of Fuzzy set
theory and thereafter it has been developed by several authors
through the contribution of the  different articles on this concept
and applied on different branches of pure and applied mathematics.
The concept of fuzzy norm was introduced by Katsaras \cite{Katsaras}
in 1984 and in 1992, Felbin\cite{Felbin1} introduced the idea of
fuzzy norm on a linear space. Cheng-Moderson \cite{Shih-chuan}
introduced another idea of fuzzy norm on a linear space whose
associated metric is same as the associated metric of
Kramosil-Michalek \cite{Kramosil}. Latter on Bag and Samanta
\cite{Bag1} modified the definition of fuzzy norm of Cheng-Moderson
\cite{Shih-chuan} and established the concept of continuity and
boundednes of a function with respect to
their fuzzy norm in\cite{Bag2}. \\
The authors T. Bag and S. K. Samanta \cite{Bag1} introduced the
definition of fuzzy norm over a linear space following the
definition S. C. Cheng and J. N. Moordeson \cite{Shih-chuan} and
they  have studied finite dimensional fuzzy normed linear spaces.
Also the definition of intuitionistic fuzzy n-normed linear space
was introduced in the paper \cite{Vijayabalaji} and established a
sufficient condition for an intuitionistic fuzzy n-normed linear
space to be complete. In this paper, following the definition of
intuitionistic fuzzy n-norm \cite{Vijayabalaji} , the definition of
generalized intuitionistic fuzzy $\psi$ norm (\, in short GIF$\psi$N
\,) is defined over a linear space. There after a sufficient
condition is given for a generalized intuitionistic fuzzy $\psi$
normed linear space to be complete and also it is proved that a
finite dimensional generalized intuitionistic fuzzy $\psi$ norm
linear space is complete. In such spaces, it is established that a
necessary and sufficient condition for a subset to be compact.
Thereafter the definition of generalized intuitionistic fuzzy $\psi$
continuity, strongly intuitionistic fuzzy $\psi$ continuity and
sequentially intuitionistic fuzzy $\psi$ continuity are defined and
proved that the concept of intuitionistic fuzzy $\psi$ continuity
and
 sequentially intuitionistic fuzzy $\psi$ continuity are
equivalent. There after it is shown that intuitionistic fuzzy
continuous image of a compact set is again a compact set.


\section{Preliminaries}

We quote some definitions and statements of a few theorems which will be
needed in the sequel.

\begin{Definition} \cite{Schweizer}.
A binary operation \, $\ast \; : \; [\,0 \; , \; 1\,] \; \times \;
[\,0 \; , \; 1\,] \;\, \longrightarrow \;\, [\,0 \; , \; 1\,]$ \, is
continuous \, $t$ - norm if \,$\ast$\, satisfies the
following conditions \, $:$ \\
$(\,i\,)$ \hspace{0.5cm} $\ast$ \, is commutative and associative ,
\\ $(\,ii\,)$ \hspace{0.4cm} $\ast$ \, is continuous , \\
$(\,iii\,)$ \hspace{0.2cm} $a \;\ast\;1 \;\,=\;\, a \hspace{1.2cm}
\forall \;\; a \;\; \varepsilon \;\; [\,0 \;,\; 1\,]$ , \\
$(\,iv\,)$ \hspace{0.2cm} $a \;\ast\; b \;\, \leq \;\, c \;\ast\; d$
\, whenever \, $a \;\leq\; c$  ,  $b \;\leq\; d$  and  $a \, , \, b
\, , \, c \, , \, d \;\, \varepsilon \;\;[\,0 \;,\; 1\,]$.
\end{Definition}

\begin{Definition}
\cite{Schweizer}. A binary operation \, $\diamond \; : \; [\,0 \; ,
\; 1\,] \; \times \; [\,0 \; , \; 1\,] \;\, \longrightarrow \;\,
[\,0 \; , \; 1\,]$ \, is continuous \, $t$-conorm if \,$\diamond$\,
satisfies the
following conditions \, $:$ \\
$(\,i\,)\;\;$ \hspace{0.1cm} $\diamond$ \, is commutative and
associative ,
\\ $(\,ii\,)\;$ \hspace{0.1cm} $\diamond$ \, is continuous , \\
$(\,iii\,)$ \hspace{0.1cm} $a \;\diamond\;0 \;\,=\;\, a
\hspace{1.2cm}
\forall \;\; a \;\; \in\;\; [\,0 \;,\; 1\,]$ , \\
$(\,iv\,)$ \hspace{0.1cm} $a \;\diamond\; b \;\, \leq \;\, c
\;\diamond\; d$ \, whenever \, $a \;\leq\; c$  ,  $b \;\leq\; d$
 and  $a \, , \, b \, , \, c \, , \, d \;\; \in\;\;[\,0
\;,\; 1\,]$.
\end{Definition}

\begin{Corollary}
\cite{Vijayabalaji}. $(\,a\,)\;$  For any \, $r_{\,1} \; , \;
r_{\,2} \;\; \in\;\; (\,0 \;,\; 1\,)$ \, with \, $r_{\,1} \;>\;
r_{\,2}$, there exist $r_{\,3} \; , \; r_{\,4} \;\; \in \;\; (\,0
\;,\; 1\,)$ \, such that \, $r_{\,1} \;\ast\; r_{\;3} \;>\; r_{\,2}$
\, and \, $r_{\,1} \;>\; r_{\,4} \;\diamond\; r_{\,2}.$
\\ $(\,b\,)$ \, For any \, $r_{\,5} \;\,
\in\;\, (\,0 \;,\; 1\,)$ , there exist \, $r_{\,6} \; , \; r_{\,7}
\;\, \in\;\, (\,0 \;,\; 1\,)$ \, such that \, $r_{\,6} \;\ast\;
r_{\,6} \;\geq\; r_{\,5}$ \,and\, $r_{\,7} \;\diamond\; r_{\,7}
\;\leq\; r_{\,5}.$
\end{Corollary}

\begin{Definition}\,\cite{Huang}\,
By an operation\, $\circ$ \,on $\mathbf{R^{\,+}}$\, we mean a two place function
$\circ \;:\; [\,0 \;,\; \infty\,) \;\times\; [\,0 \;,\; \infty\,)\;\longrightarrow
\; [\,0 \;,\; \infty\,) $\, which is associative, commutative, non decreasing in each 
place and such that $a \;\circ\; 0 \,=\, a \;\; \forall \;a\;\in\,[\,0 \;,\; \infty\,)$. The most used operations on $\mathbf{R^{\,+}}$\, are \,$\circ_{\,a}\,(\,s \;,\; t\,)
\;=\; s \;+\; t$ \,,\, $\circ_{\,m}\,(\,s \;,\; t\,)
\;=\; \max\{s \;,\; t\}$,  $\circ_{\,n}\,(\,s \;,\; t\,)
\;=\; \left(\,s^{\,n} \;+\; t^{\,n}\,\right)^{\,\frac{1}{n}}$.

\end{Definition}

\begin{Definition}
Let $\psi$ be a function defined on the real field $\mathbf{R}$ into itself
satisfying the following properties $\textbf{:}$\\
$(\,i\,) \hspace{0.5cm} \psi(\,-\,t\,)\;=\;\psi(\,t\,)$ \,for all\, $t\,\in
\,\mathbf{R}$ \\ $(\,ii\,) \hspace{0.4cm} \psi(\,1\,)\;=\;1$\\
$(\,iii\,) \hspace{0.3cm} \psi$ is strictly increasing and continuous on
$(\,0 \;,\;\infty\,)$ \\ $(\,iv\,) \hspace{0.3cm}
\mathop {\lim }\limits_{\alpha \; \to \;0} \,\psi (\,\alpha \,)\; = \;0$
and $ \mathop {\lim }\limits_{\alpha \; \to \;\infty} \,\psi (\,\alpha \,)\;
 = \;\infty $
\end{Definition}

\begin{Example}
As example of such functions, consider\, $\psi(\,\alpha\,) \,=\,
\left| {\,\alpha \,} \right|$ ;\\ $\psi(\,\alpha\,) \,=\,
\left| {\,\alpha \,} \right|^{\,p}$ , $p\,\in\,\mathbf{R^{\,+}}$ ;
$\psi(\,\alpha\,) \,=\,\frac{2\,\alpha^{\,2\,n}}{\left|
{\,\alpha \,} \right|\;+\;1}$ , $n\,\in\,\mathbf{N^{\,+}}$. The function
$\psi$ allows us to generalize fuzzy metric and normed space.
\end{Example}

\begin{Definition}
\cite{Samanta}. Let \,$\ast$\, be a continuous \,$t$-norm ,
\,$\diamond$\, be a continuous \,$t$-conorm  and \,$V$\, be a
linear space over the field \, $F \,(\, = \, \mathbf{R}\,$ or
$\, \mathbf{C} \;)$. An \textbf{intuitionistic fuzzy norm} on
\,$V$\, is an object of the form \, $A \;\,=\;\, \{\; (\,(\,x \;,\;
t\,) \;,\; \mu\,(\,x \;,\; t\,) \;,\; \nu\,(\,x \;,\; t\,) \;) \;\,
: \;\, (\,x \;,\; t\,) \;\,\in\;\, V \;\times\; \mathbf{R^{\,+}}
\;\}$ , where $\mu \,,\, \nu\;$ are fuzzy sets on \, $V \;\times\;
\mathbf{R^{\,+}}$ , \,$\mu$\, denotes the degree of membership and
\,$\nu$\, denotes the degree of non - membership \,$(\,x \;,\; t\,)
\;\,\in\;\, V \;\times\;
\mathbf{R^{\,+}}$\, satisfying the following conditions $:$ \\
$(\,i\,)$ \hspace{0.10cm}  $\mu\,(\,x \;,\; t\,) \;+\; \nu\,(\,x
\;,\; t\,) \;\,\leq\;\, 1 \hspace{1.2cm} \forall \;\; (\,x \;,\;
t\,)
\;\,\in\;\, V \;\times\; \mathbf{R^{\,+}}\, ;$ \\
$(\,ii\,)$ \hspace{0.10cm}$\mu\,(\,x \;,\; t\,) \;\,>\;\, 0 \, ;$ \\
$(\,iii\,)$ $\mu\,(\,x \;,\; t\,) \;\,=\;\, 1$ \, if
and only if \, $x \;=\; \theta \,$, $\theta$ is null vector ; \\
$(\,iv\,)$\hspace{0.05cm} $\mu\,(\,c\,x \;,\; t\,) \;\,=\;\,
\mu\,(\,x \;,\; \frac{t}{|\,c\,|}\,)$ \, \, $\;\forall\; c
\;\,\in\;\, F \, $ and $c \;\neq\; 0 \;;$ \\ $(\,v\,)$
\hspace{0.10cm} $\mu\,(\,x \;,\; s\,) \;\ast\; \mu\,(\,y \;,\; t\,)
\;\,\leq\;\, \mu\,(\,x \;+\; y \;,\; s \;+\; t\,) \, ;$ \\
$(\,vi\,)$ \hspace{0.05cm} $\mu\,(\,x \;,\; \cdot\,)$ is
non-decreasing function of \, $\mathbf{R^{\,+}}$ \,and\,
$\mathop{\lim }\limits_{t\;\, \to \,\;\infty } \;\mu\,\left(
{\;x\;,\;t\,} \right)=1 ;$
\\ $(\,vii\,)$ \hspace{0.10cm}$\nu\,(\,x \;,\; t\,) \;\,<\;\, 1 \, ;$ \\
$(\,viii\,)$ $\nu\,(\,x \;,\; t\,) \;\,=\;\, 0$ \, if and only if \,
$x \;=\; \theta \, ;$ \\ $(\,ix\,)$ \hspace{0.05cm} $\nu\,(\,c\,x
\;,\; t\,) \;\,=\;\, \nu\,(\,x \;,\; \frac{t}{|\,c\,|}\,)$ \, \,
$\;\forall\; c \;\,\in\;\, F \, $ and $c \;\neq\; 0 \;;$ \\
$(\,x\,)$ \hspace{0.15cm} $\nu\,(\,x \;,\; s\,) \;\diamond\;
\nu\,(\,y \;,\; t\,) \;\,\geq\;\, \nu\,(\,x \;+\; y \;,\; s \;+\;
t\,) \, ;$ \\ $(\,xi\,)$ \hspace{0.04cm} $\nu\,(\,x \;,\; \cdot\,)$
s non-increasing function of \, $\mathbf{R^{\,+}}$ \,and\, $\mathop
{\lim }\limits_{t\;\, \to \,\;\infty } \;\,\,\nu\,\left(
{\;x\;,\;t\,} \right)=0.$
\end{Definition}

\smallskip
\begin{Definition}
\cite{Samanta}. If $A$ is an intuitionistic fuzzy norm on a linear
space $V$ then $(V\;,\;A)$ is called an intuitionistic fuzzy normed
linear space.
\end{Definition}

\smallskip
For the intuitionistic fuzzy normed linear space \,$(\,V \;,\;
A\,)$\,, we further assume that \,$\mu,\, \nu,\, \ast,\, \diamond$\, satisfy
the following axioms : \newline
$(\,xii\,)$\hspace{0.53cm} $\left. {{}_{a\;\; \ast \;\;a\;\; = \;\;a}^{a\;\;
\diamond \;\;a\;\; = \;\;a} \;\;} \right\}\;\;\;$\hspace{0.5cm},\,\,\,for
all $\;a\;\; \in \;\;[\,0\;\,,\;\,1\,].$ \newline
$(\,xiii\,)$ \,\, $\mu\,(\,x \;,\; t\,) \;>\; 0 \;\;\;\;,$\, for all $\;\; t
\;>\; 0 \;\; \Rightarrow \;\; x \;=\;\theta\;.$ \newline
$(\,xiv\,)$\,\,\, $\nu\,(\,x \;,\; t\,) \;<\; 1 \;\;\;\;\;\; ,$ \,for all $%
\;\; t \;>\; 0 \;\; \Rightarrow \;\; x \;=\; \theta\;.$\newline

\begin{Definition}
\cite{Samanta}. A sequence $\{\,x_n\,\}_n$ in an intuitionistic
fuzzy normed linear space $(\,V\,,\,A\,)$ is said to
\textbf{converge} to $x\;\in\;V$ if for given $r>0,\;t>0,\;0<r<1$,
there exist an integer
$n_0\;\in\;\mathbf{N}$ such that \\
$\;\mu\,(\,x_n\,-\,x\,,\,t\,)\;>\;1\,-\,r$
 \,\,and\,\, $\nu\,(\,x_n\,-\,x\,,\,t\,)\;<\;r$ \,\,for all $n\;\geq \;n_0$.
\end{Definition}
\smallskip

\begin{Definition}
\cite{Samanta}. A sequence $\{\,x_n\,\}_n$ in an intuitionistic
fuzzy normed linear space $(\,V\,,\,A\,)$ is said to be
\textbf{cauchy sequence} if $\mathop {\lim }\limits_{n\;\, \to
\,\;\infty } \;\mu\,(x_{n+p}-x_n ,t)\;=\;1\; $ and $\mathop {\lim
}\limits_{n\;\, \to \,\;\infty } \;\nu\,(x_{n+p}-x_n
,t)\;=\;0\;\;,\;p\;=\;1\,,\,2\,,\,3\,,\;\;\cdots $.
\end{Definition}
\smallskip

\begin{Definition}
\cite{Samanta}. Let, $(\,U\;,\;A\,)$ and $(\,V\;,\;B\,)$ be two
intuitionistic fuzzy normed linear space over the same field $F$. A
mapping $f$ from $(\,U\;,\;A\,)$ to $(\,V\;,\;B\,)$ is said to be
\textbf{ intuitionistic fuzzy continuous} at $x_0\;\in\;U$, if for
any given
$\epsilon\;>\;0\;,\;\alpha\;\in\;(\,0\;,\;1\,)\;,\;\;\exists\;\delta
\;=\;\delta(\,\alpha\,,\,\epsilon\,)\;>0\;,\;\beta\;=\;
\beta(\,\alpha\,,\,\epsilon\,)\;\in\;(\,0\,,\,1\,)$ \,such that for
all $x\;\in\;U$,
\[\mu_U(x-x_0 \;,\;\delta)\;>\;1-\beta\;\Rightarrow\;
\mu_V(f(x)-f(x_0) \;,\;\epsilon)\;>\;1-\alpha\] \[ \nu_U(x-x_0
\;,\;\delta)\;<\;\beta\;\Rightarrow\; \nu_V(f(x)-f(x_0)
\;,\;\epsilon)\;<\;\alpha\;. \]
\end{Definition}
\smallskip

\begin{Definition}
\cite{Samanta}. A mapping $f$ from $(\,U\,,\,A\,)$ to
$(\,V\,,\,B\,)$ is said to be \textbf{sequentially intuitionistic
fuzzy continuous} at $x_0\;\in\;U$, if for any sequence
$\{\,x_n\,\}_n$, $x_n\;\in\;U\;,\;\forall\;n\;\in\;\mathbf{N}$ with
$x_n\;\rightarrow\;x_0$ in $(\,U\,,\,A\,)$ implies
$f(x_n)\;\rightarrow\;f(x_0)$ in $(\,V\,,\,B\,)$, that is
\[\mathop {\lim }\limits_{n\;\,
\to \,\;\infty } \;\mu_U(x_n-x_0 \,,\,t)\;=\;1 \;and \; \mathop
{\lim }\limits_{n\;\, \to \,\;\infty } \;\nu_U(x_n-x_0
\,,\,t)\;=\;0\; \] \[\Rightarrow\;\mathop {\lim }\limits_{n\;\, \to
\,\;\infty } \;\mu_{\,V}(f(x_n)-f(x_0)\,,\,t)\;=\;1 \; and \;
\mathop {\lim }\limits_{n\;\, \to \,\;\infty }
\;\nu_{\,V}(f(x_n)-f(x_0)\,,\,t)\;=\;0\;.\]
\end{Definition}
\smallskip

\begin{Theorem}
\cite{Samanta}.  Let  $f$  be a mapping from $(\,U\,,\,A\,)$ to
$(\,V\,,\,B\,)$. Then $f$ is intuitionistic fuzzy continuous on $U$
if and only if it is sequentially intuitionistic fuzzy continuous on
$U$.
\end{Theorem}
\bigskip


\section{ Generalized Intuitionistic Fuzzy \\ $\psi$-Normed Linear space}

\begin{Definition}
 Let \,$\ast$\, be a continuous \,$t$-norm ,
\,$\diamond$\, be a continuous \,$t$-conorm  and \,$V$\, be a
linear space over the field \, $\mathbf{R}$.
 A \textbf{Generalized intuitionistic fuzzy
$\psi$-norm} on \,$V$\, is an object of the form \,
${\;A^{\,\psi } \,}\;=\; \{\; (\,(\,x \;,\;
t\,) \;,\\ \mu\,(\,x \;,\; t\,) \;,\; \nu\,(\,x \;,\; t\,) \;) \;\,
: \;\, (\,x \;,\; t\,) \;\,\in\;\, V \;\times\; \mathbf{R^{\,+}}
\;\}$ , where $\mu \,,\, \nu\;$ are fuzzy sets on \, $V \;\times\;
\mathbf{R^{\,+}}$ , \,$\mu$\, denotes the degree of membership and
\,$\nu$\, denotes the degree of non-membership \,$(\,x \;,\; t\,)
\;\,\in\; V \;\times\;
\mathbf{R^{\,+}}$\, satisfying the following conditions $:$ \\
$(\,i\,)$ \hspace{0.10cm}  $\mu\,(\,x \;,\; t\,) \;+\; \nu\,(\,x
\;,\; t\,) \;\,\leq\;\, 1 \hspace{1.2cm} \forall \;\; (\,x \;,\;
t\,)
\;\,\in\;\, V \;\times\; \mathbf{R^{\,+}}\, ;$ \\
$(\,ii\,)$ \hspace{0.10cm}$\mu\,(\,x \;,\; t\,) \;\,>\;\, 0 \, ;$ \\
$(\,iii\,)$ $\mu\,(\,x \;,\; t\,) \;\,=\;\, 1$ \, if
and only if \, $x \;=\; \theta \,$, $\theta$ is null vector ; \\
$(\,iv\,)$\hspace{0.05cm} $\mu\,(\,\alpha\,x \;,\; t\,) \;\,=\;\,
\mu\,(\,x \;,\; \frac{t}{\psi(\,\alpha\,)}\,)$ \, \, $\forall\;\alpha
\;\,\in\;\, \mathbf{R} \, $ and $\alpha\;\neq\;0$ \\ $(\,v\,)$
\hspace{0.10cm} $\mu\,(\,x \;,\; s\,) \;\ast\; \mu\,(\,y \;,\; t\,)
\;\,\leq\;\, \mu\,(\,x \;+\; y \;,\; s\;\circ\;t\,
\,) \, ;$ \\
$(\,vi\,)$ \hspace{0.05cm} $\mu\,(\,x \;,\; \cdot\,)$ is
non-decreasing function of \, $\mathbf{R^{\,+}}$ \,and\,
$\mathop{\lim }\limits_{t\;\, \to \,\;\infty } \;\mu\,\left(
{\;x\;,\;t\,} \right)=1 ;$
\\ $(\,vii\,)$ \hspace{0.10cm}$\nu\,(\,x \;,\; t\,) \;\,<\;\, 1 \, ;$ \\
$(\,viii\,)$ $\nu\,(\,x \;,\; t\,) \;\,=\;\, 0$ \, if and only if \,
$x \;=\; \theta \, ;$ \\ $(\,ix\,)$ \hspace{0.05cm}
$\nu\,(\,\alpha\,x \;,\; t\,) \;\,=\;\,\nu\,(\,x \;,\;
 \frac{t}{\psi(\,\alpha\,)}\,)$ \, \, $\forall\;\alpha
\;\,\in\;\, \mathbf{R} \, $ and $\alpha\;\neq\;0$ \\
$(\,x\,)$ \hspace{0.15cm} $\nu\,(\,x \;,\; s\,) \;\diamond\;
\nu\,(\,y \;,\; t\,) \;\,\geq\;\, \nu\,(\,x \;+\; y \;,\;
 s\;\circ\;t\,) \, ;$ \\ $(\,xi\,)$
\hspace{0.04cm} $\nu\,(\,x \;,\; \cdot\,)$ s non-increasing function of \, $\mathbf{R^{\,+}}$ \,and\, $\mathop
{\lim }\limits_{t\;\, \to \,\;\infty } \;\,\,\nu\,\left(
{\;x\;,\;t\,} \right)=0.$ \\\\
If $A$ is an Generalized intuitionistic fuzzy
$\psi$-norm on a linear
space $V$ then $(V\;,\;A)$ is called a Generalized intuitionistic fuzzy
$\psi$-normed linear space.
\end{Definition}

\begin{Definition}
A sequence $\{\,x_n\,\}_n$ in a generalized IF$\psi$NLS
$\left( {\,V\;,\;A^{\,\psi } \,} \right)$ is said to converge to $x\in\,V$ if for any given $r>0$ , $t>0$ , $r\,\in\,(\,0\,,\,1\,)$ there exists an integer $n_0\,\in\,N$ such that $\;\mu\,(\,x_n\,-\,x\,,\,t\,)\;>\;1\,-\,r$\,\,and\,\, $\nu\,(\,x_n\,-\,x\,,\,t\,)\;<\;r\;\;\;\forall\;n\,
\geq\,n_0\;\,$.
\end{Definition}

\begin{Theorem}
In a Generalized intuitionistic fuzzy
$\psi$-normed linear space $(V\;,\;A^{\,\psi })$, a sequence $\{\,x_{\,n}\}$ converges to $x$
if and only if $\mu(\,x_{\,n} \,-\, x \;,\; t\,) \;\longrightarrow \;1$\, and \,$\nu(\,x_{\,n} \,-\, x \;,\; t\,) \;\longrightarrow \;0$\, as $ n \;\rightarrow\;\infty$.
\end{Theorem}
{\bf Proof.}$\;\;$
Fix $t\,>\,0\;$. Suppose $\{\,x_{\,n}\}$ converges to $x$ in $(V\;,\;A^{\,\psi })$ . Then for a given $r,\;r\,\in\,(\,0\,,\,1\,)$, there exists an integer $n_0\,\in\,N$ such that $\;\mu\,(\,x_n\,-\,x\,,\,t\,)\;>\;1\,-\,r$\,\,and\,\, $\nu\,(\,x_n\,-\,x\,,\,t\,)\;<\;r$ \,. Thus $\;1\,-\,\mu\,(\,x_n\,-\,x\,,\,t\,)\;<\,r$\,\,\,and $\nu\,(\,x_n\,-\,x\,,\,t\,)\;<\;r$ \,, and hence $\mu(\,x_{\,n} \,-\, x \;,\; t\,) \;\longrightarrow \;1$\, and \,$\nu(\,x_{\,n} \,-\, x \;,\; t\,) \;\longrightarrow \;0$ \,as \,$ n \;\rightarrow\;\infty$.
Conversely , if for each $t\,>\,0\;$ , $\mu(\,x_{\,n} \,-\, x \;,\; t\,) \;\longrightarrow \;1$\, and \,$\nu(\,x_{\,n} \,-\, x \;,\; t\,) \;\longrightarrow \;0$\, as  $ n \;\rightarrow\;\infty$, then for every $r,\;r\,\in\,(\,0\,,\,1\,)$, there exists an integer $n_0\,$ such that $\;1\,-\,\mu\,(\,x_n\,-\,x\,,\,t\,)\;<\,r$\,\,\,and $\nu\,(\,x_n\,-\,x\,,\,t\,)\;<\;r$ \,$\forall\;n\,
\geq\,n_0\;\,$. Thus $\;\mu\,(\,x_n\,-\,x\,,\,t\,)\;>\;1\,-\,r$
 \,\,and\,\, $\nu\,(\,x_n\,-\,x\,,\,t\,)\;<\;r$ \,\,for all $n\;\geq \;n_0$. Hence $\{\,x_{\,n}\}$ converges to $x$ in $(V\;,\;A^{\,\psi })$.

\medskip

\begin{Theorem}
The limit is unique for a convergent sequence $\{\,x_n\,\}_n$ in a generalized
IF$\psi$NLS $\left( {\,V\;,\;A^{\,\psi } \,} \right)$ .
\end{Theorem}

{\bf Proof.}$\;\;$
Let $\mathop {\lim }\limits_{n\;\, \to \;\,\infty } \;x_{\,n} \;\; =
\;\;x$ \,and \,$\mathop {\lim }\limits_{n\;\, \to \;\,\infty }
\;x_{\,n} \;\; = \;\;y$ .
Also let $
 s\,,\,t \;\in \mathbf{R^{\,+}}$. Now,\[\mathop {\lim }\limits_{n\;\,
\to \;\,\infty } \;x_{\,n} \;\; = \;\;x\;\; \Rightarrow \;\;\left\{
{_{\mathop {\lim }\limits_{n\;\, \to \;\,\infty } \;\nu\,\left(
{\,x_{\,n} \; - \;x\;,\;t\,} \right)\;\; = \;\;0}^{\mathop {\lim
}\limits_{n\;\, \to \;\,\infty } \;\mu\,\left( {\,x_{\,n} \; -
\;x\;,\;t\,} \right)\;\; = \;\;1} } \right.\]  \[\mathop {\lim
}\limits_{n\;\, \to \;\,\infty } \;x_{\,n} \;\; = \;\;y\;\;
\Rightarrow \;\;\left\{ {_{\mathop {\lim }\limits_{n\;\, \to
\;\,\infty } \;\nu\,\left( {\,x_{\,n} \; - \;y\;,\;t\,} \right)\;\; =
\;\;0}^{\mathop {\lim }\limits_{n\;\, \to \;\,\infty } \;\mu\,\left(
{\,x_{\,n} \; - \;y\;,\;t\,} \right)\;\; = \;\;1} }\right.\] \\
$\begin{array}{l}
 \nu\,\left( {\,x\; - \;y\;,
 s\;\circ\;t\,} \right)\;\; = \;\;\nu\,\left(
 {\,x\; - \;x_{\,n} \; + \;x_{\,n} \; - \;y\;,\;,
 s\;\circ\;t\,} \right) \\
 {\hspace{3.9cm}} \le \;\;\nu\,\left({\,x\; -
 \;x_{\,n} \;,\;s\,} \right)\;\diamond\;\nu\,\left( {\,x_{\,n} \; -
 \;y\;,\;t\,} \right) \\
{\hspace{3.9cm} = }\;\;\nu\,
 \left(-({\,x_{\,n} \; - \;x\;),\;s\,} \right)\;\diamond\;\nu\,\left( {\,x_{\,n} \; - \;y\;,\;t\,} \right) \\
 {\hspace{3.9cm} = }\;\;\nu\,
 \left( {\,x_{\,n} \; - \;x\;,\frac{s}{\psi(\,\;-1\,)}\,} \right)\;\diamond\;\nu\,\left( {\,x_{\,n} \; - \;y\;,\;t\,} \right) \\
 {\hspace{3.9cm} = }\;\;\nu\,
 \left( {\,x_{\,n} \; - \;x\;,\frac{s}{\psi(\,\-1\,)}\,} \right)\;\diamond\;\nu\,\left( {\,x_{\,n} \; - \;y\;,\;t\,} \right) \\
 {\hspace{3.9cm} = }\;\;\nu\,
 \left( {\,x_{\,n} \; - \;x\;,\;s\,} \right)\;\diamond\;\nu\,\left( {\,x_{\,n} \; - \;y\;,\;t\,} \right) \\
 \end{array}$
\\Taking limit , we have \\
\[\nu\,(\,x\; - \;y\;,
 s\;\circ\;t\,) \;\,\leq\;\,
\mathop {\lim }\limits_{n\; \to \;\infty } \;\nu\,\left( {\,x_{\,n} \;
- \;x\;,\;s\,} \right) \;\diamond\; \mathop {\lim }\limits_{n\; \to
\;\infty } \;\nu\,\left( {\,x_{\,n} \; - \;y\;,\;t\,} \right)
\;\,=\;\, 0\] \[\Longrightarrow \;\; \nu\,(\,x\; - \;y\;,
 s\;\circ\;t\,) \;\,=\;\, 0\;\;\Longrightarrow\;\; x \;-\; y \;\,=\;\,
\underline{0} \;\;\Longrightarrow\;\; x \;\,=\;\,y\]

\begin{Theorem}
If $\mathop {\lim }\limits_{n\;\, \to \;\,\infty } \;x_{\,n} \;\; =
\;\;x$\,\,and \,$\mathop {\lim }\limits_{n\;\, \to \;\,\infty }
\;y_{\,n} \;\;=\;\;y$ \,then\, $\mathop {\lim }\limits_{n\;\, \to \;\,\infty } \;\left(\,x_{\,n} \;\;  +
\;y_{\,n}\,\right) \;\;=\;\;x+y$ \,in a generalized IF$\psi$NLS
\end{Theorem}

{\bf Proof.}$\;$The proof directly follows from from the proof of the theorem 3 \cite{Samanta}

\begin{Theorem}
\,$\mathop {\lim }\limits_{n\;\, \to \;\,\infty } \;x_{\,n} \;\;
= \;\;x$\, and \,$k \;(\,\neq\; 0\,) \;\,\varepsilon \;\, F$\, $\Longrightarrow$
\,$\mathop {\lim }\limits_{n\;\, \to \;\,\infty } \;k\,x_{\,n} \;\;
= \;\;k\,x$\,\, in a generalized IF$\psi$NLS
$\left( {\,V\;,\;A^{\,\psi } \,} \right)$
\end{Theorem}

{\bf Proof.}$\;\;$Obvious.

\begin{Definition}
\ A sequence $\{\,x_n\,\}_n$ in a generalized IF$\psi$NLS
$\left( {\,V\;,\;A^{\,\psi } \,} \right)$  is said to be
\textbf{cauchy sequence} if for any given $r>0$, $t>0$,
$\;r\,\in\,(\,0\,,\,1\,)$ there exists an integer $n_0\,\in\,N$ such that $\;\mu\,(\,x_m\,-\,x_n\,,\,t\,)\;>\;1\,-\,r$\,\,and\,\, $\nu\,(\,x_m\,-\,x_n\,,\,t\,)\;<\;r\;\;\;\forall\;m \,,\,n\;
\geq\,n_0\;\,$.
\end{Definition}

\begin{Theorem}
In a generalized IF$\psi$NLS
$\left( {\,V\;,\;A^{\,\psi } \,} \right)$ , every convergent sequence is a
Cauchy sequence.
\end{Theorem}

{\bf Proof.}$\;\;$
Let \,$\left\{ {\;x_{\,n}\;} \right\}_{\,n} $\, be a convergent
sequence in the \,IF$\psi$NLS
$\left( {\,V\;,\;A^{\,\psi } \,} \right)$, with \,$\mathop {\lim
}\limits_{n\;\, \to \;\,\infty } \;x_{\,n} \;\; = \;\;x$.\\Let $
 (\,
 s\;\circ\;t\,)\in \mathbf{R^{\,+}}$\, and
 $\;p\;=\;1\,,\,2\,,\,3\,,\;\;\cdots $.\\we have \\$\begin{array}{l}
 \mu\,\left( {\,x_{\,n\; + \;p} \; - \;x_{\,n} \;,\;s\;\circ\;t\,} \right)
 \;\; = \;\;\mu\,\left( {\,x_{\,n\; + \;p} \; - \;x\; + \;x\; - \;x_{\,n}
  \;,\;s\;\circ\;t\,} \right) \\
 {\hspace{5.0cm}} \ge \;\;\mu\,\left( {\,x_{\,n\; + \;p} \; - \;x\;,\;s\,}
  \right)\;\, * \;\,\mu\,\left({\,x\; -
 \;x_{\,n} \;,\;t\,} \right)\;
 {\hspace{5.1cm} =}\;\;\mu\,\left( {\,x_{\,n\; + \;p} \; - \;x\;,\;s\,} \right)
 \; * \;\mu\,\left(-({\,x_{\,n} \; - \;x\;),\;s\,} \right)\; \\
 \end{array}$ \\$\begin{array}{l}
  {\hspace{5.0cm}=}\;\;\mu\,\left( {\,x_{\,n\; + \;p} \; - \;x\;,\;s\,}
  \right)\;\, * \;\,\mu\,\left(-({\,x_{\,n} \; - \;x\;),\;t\,} \right)\;
 {\hspace{5.1cm} =}\;\;\mu\,\left( {\,x_{\,n\; + \;p} \; - \;x\;,\;s\,} \right)
 \; * \;\mu\,\left( {\,x_{\,n} \; - \;x\;,\;t\,} \right) \\{\hspace{5.1cm} =}\;\;\mu\,\left( {\,x_{\,n\; + \;p} \; - \;x\;,\;s\,} \right)
 \; * \;\mu\,\left( {\,x_{\,n} \; - \;x\;,\frac{t}{\psi(\,\;-1\,)}\,} \right)\;\\{\hspace{5.1cm} =}\;\;\mu\,\left( {\,x_{\,n\; + \;p} \; - \;x\;,\;s\,} \right)
 \; * \;\mu\,\left( {\,x_{\,n} \; - \;x\;,\frac{t}{\psi(\,\-1\,)}\,} \right)\; \\{\hspace{5.1cm} =}\;\;\mu\,\left( {\,x_{\,n\; + \;p} \; - \;x\;,\;s\,} \right)
 \; *\;\mu\,\left( {\,x_{\,n} \; - \;x\;,\;t\,}\right)\; \\
 \end{array}$ \\\\ Let $r > 0$, $t \, , \, s > 0$,
$\,r\,\in\,(\,0\,,\,1\,)$, then $\exists$ an integer $n_0\,\in\,N$ such that
 \\\[ \begin{array}{l}
 \mu\,\left( {\,(\,x_{\,n\; + \;p} \; - \;x_{\,n}\,)\;\,,\,\;s\;\circ\;t\,}
 \right) \\ {\hspace{3.5cm}} \ge \;\;
  \mu\,\left( {\,x_{\,n\,+\,p} \; - \;x\;,\;s\,} \right)\;\; *
 \;\;\mu\,\left( {\,x_{\,n}
 \; - \;x\;,\;t\,} \right) \\
 {\hspace{3.5cm}} = \;\;(\,1\;-\;r\,)\; * \;(\,1\;-\;r\,)\;\; = \;\;(\,1\;-\;r\,)
 \;\; \forall\;n\,\geq\,n_0\\
 \end{array}
\]
 Again,\\ $\begin{array}{l}
 \nu\,\left( {\,x_{\,n\; + \;p} \; - \;x_{\,n} \;,\;s\; \circ \;t\,} \right)
 \;\; = \;\;\nu\,\left( {\,x_{\,n\; + \;p} \; - \;x\; + \;x\; - \;x_{\,n}
  \;,\;s\; \circ \;t\,} \right) \\
 {\hspace{5.0cm}} \leq\;\;\nu\,\left( {\,x_{\,n\; + \;p} \; - \;x\;,\;s\,}
  \right)\;\, \diamond\ \;\,\nu\,\left({\,x\; -
 \;x_{\,n} \;,\;t\,} \right)\;
 {\hspace{5.1cm} =}\;\;\nu\,\left( {\,x_{\,n\; + \;p} \; - \;x\;,\;s\,} \right)
 \; * \;\nu\,\left(-({\,x_{\,n} \; - \;x\;),\;s\,} \right)\; \\
 \end{array}$\\$\begin{array}{l}
  {\hspace{5.0cm}=}\;\;\nu\,\left( {\,x_{\,n\; + \;p} \; - \;x\;,\;s\,}
  \right)\;\, \diamond\ \;\,\nu\,\left(-({\,x_{\,n} \; - \;x\;),\;t\,} \right)\;
 {\hspace{5.1cm} =}\;\;\nu\,\left( {\,x_{\,n\; + \;p} \; - \;x\;,\;s\,} \right)
 \; \diamond\ \;\nu\,\left( {\,x_{\,n} \; - \;x\;,\;t\,} \right) \\{\hspace{5.1cm} =}\;\;\nu\,\left( {\,x_{\,n\; + \;p} \; - \;x\;,\;s\,} \right)
 \; \diamond\ \;\nu\,\left( {\,x_{\,n} \; - \;x\;,\frac{t}{\psi(\,\;-1\,)}\,} \right)\;\\{\hspace{5.1cm} =}\;\;\nu\,\left( {\,x_{\,n\; + \;p} \; - \;x\;,\;s\,} \right)
 \; \diamond\ \;\nu\,\left( {\,x_{\,n} \; - \;x\;,\frac{t}{\psi(\,\-1\,)}\,} \right)\; \\{\hspace{5.1cm} =}\;\;\nu\,\left( {\,x_{\,n\; + \;p} \; - \;x\;,\;s\,} \right)
 \; \diamond\ \;\nu\,\left( {\,x_{\,n} \; - \;x\;,\;t\,}\right)\; \\
 \end{array}$ \\\\ Thus, we see that \\\[
\begin{array}{l}
 \nu\,\left( {\,(\,x_{\,n\; + \;p} \; - \;x_{\,n}\,)\;\,,\,\;s\;\circ\;t\,}
 \right) \\ {\hspace{3.5cm}} \leq\;\; \nu\,
 \left( {\,x_{\,n\,+\,p} \; - \;x\;,\;s\,} \right)\;\; \diamond\
 \;\;\nu\,\left( {\,x_{\,n}
 \; - \;x\;,\;t\,} \right) \\
 {\hspace{3.5cm}} = \;\;r\; \diamond\ \;r\;\; = \;\;r \;\; \forall\;n\,\geq\,n_0 \\
 \end{array}
\]\\
Thus, \,$\left\{ {\;x_{\,n} \;} \right\}_{\,n} $\, is a Cauchy
sequence In a Generalized IF$\psi$NLS $\left( {\,V\;,\;A^{\,\psi } \,} \right)$.

\begin{Note}
The converse of the above theorem is not necessarily true . It can be
verified by the following example .
\end{Note}

\begin{Example}
Let \,$(\,V \;,\; \|\,\cdot\,\| \,)$ \, be a normed linear space and
define \,$a \;\ast\; b \;\,=\;\, \min\,\{\;a \;,\; b\;\}$ \,and  \,$a
\;\diamond\; b \;\,=\;\, \max\,\{\;a \;,\; b\;\}$ for all $\;a,b\,\in\,(\,0\,,\,1\,)$.
For all \,$t \;>\;
0$, define \, $\mu\,(\,x \;,\; t\,) \;\,=\;\, \frac{t}{t \;+\;
k\;\|\,x\,\|}$\,,\, $\nu\,(\,x \;,\; t\,) \;\,=\;\,
\frac{k\;\|\,x\,\|}{t \;+\; k\;\|\,x\,\|}$\, where \, $k \;>\; 0$ and
$\psi(\,t\,) \,=\,|\,t\,|$ .
It is easy to see that ${\,A^{\,\psi } \,} \;=\; \{\; (\,(\,x \;,\;
t\,) \;, \mu\,(\,x \;,\; t\,) \;,\; \nu\,(\,x \;,\; t\,) \;) \;\,
: \;\, (\,x \;,\; t\,) \;\,\in\;\, V \;\times\; \mathbf{R^{\,+}}
\;\}$ is a generalized IF$\psi$NLS.  Then \\ $(\,i\,)$ \,$\left\{ {\;x_{\,n}\;} \right\}_{\,n} $\,is a Cauchy sequence in \,$(\,V
\;,\; \|\,\cdot\,\| \,)$ \, if and only if \,$\left\{ {\;x_{\,n}\;} \right\}_{\,n} $\,is a Cauchy sequence in a Generalized IF$\psi$NLS $\left( {\,V\;,\;A^{\,\psi } \,} \right)$. \\ $(\,ii\,)$ \,$\left\{ {\;x_{\,n}\;} \right\}_{\,n} $\,is a convergent sequence in  \,$(\,V
\;,\; \|\,\cdot\,\| \,)$ \, if and only if \,$\left\{ {\;x_{\,n}\;} \right\}_{\,n} $\,
is a convergent sequence in a Generalized IF$\psi$NLS $\left( {\,V\;,\;A^{\,\psi } \,} \right)$.
\end{Example}

{\bf Proof.}$\;\;$The verification directly follows from the example 2 of \cite{Samanta}.

\begin{Theorem}
In a generalized IF$\psi$NLS
$\left( {\,V\;,\;A^{\,\psi } \,} \right)$ , a sequence  \,$\left\{ {\;x_{\,n} \;} \right\}_{\,n} $\, is a Cauchy
sequence if and only if  $\mu\,(\,x_{n\,+\,p}\,-\,x_n \,,\, t\,) \;\longrightarrow \;1$\, and \,$\nu\,(\,x_{n\,+\,p}\,-\,x_n
\,,\, t\,) \;\longrightarrow \;0$\, as\, $ n \;\rightarrow\;\infty$.
\end{Theorem}

{\bf Proof.}$\;\;$
Fix $t\,>\,0\;$. Suppose $\{\,x_{\,n}\}$ is a Cauchy sequence in $\left( {\,V\;,\;A^{\,\psi } \,} \right)$. Then for a given $r \,>\, 0,\;r\,\in\,(\,0\,,\,1\,)$, there exists an integer $n_0\,\in\,N$ such that $\;\mu\,(\,x_{n\,+\,p}\,-\,x_n \,,\, t\,)\;>\;1\,-\,r$\,\,and\,\, $\nu\,(\,x_{n\,+\,p}\,-\,x_n
\,,\, t\,)$. Thus $\;1\,-\,\mu\,(\,x_{n\,+\,p}\,-\,x_n \,,\, t\,)\;<\,r$\,\,\,and $\nu\,(\,x_{n\,+\,p}\,-\,x_n
\,,\, t\,)\;<\;r$ , and hence $\mu\,(\,x_{n\,+\,p}\,-\,x_n \,,\, t\,) \;\longrightarrow \;1$\, and \,$\nu\,(\,x_{n\,+\,p}\,-\,x_n
\,,\, t\,) \;\longrightarrow \;0$\,as $ n \;\rightarrow\;\infty$.\\
Conversely , for each $t\,>\,0\;$, suppose $\mu\,(\,x_{n\,+\,p}\,-\,x_n \,,\, t\,) \;\longrightarrow \;1$\, and \,$\nu\,(\,x_{n\,+\,p}\,-\,x_n
\,,\, t\,) \;\longrightarrow \;0$\, as $ n \;\rightarrow\;\infty$. Then for every $r\,>\,0 ,\;r\,\in\,(\,0\,,\,1\,)$, there exists an integer \,$n_0\,\in\,N$\, such that $\;1\,-\,\mu\,(\,x_{n\,+\,p}\,-\,x_n \,,\, t\,)\;<\;r$\,\,\,and $\nu\,(\,x_{n\,+\,p}\,-\,x_n
\,,\, t\,)\;<\;r$ \,$\forall\;n\,
\geq\,n_0$. Thus $\;\mu\,(\,x_{n\,+\,p}\,-\,x_n \,,\, t\,)\;>\;1\,-\,r$
 \,\,and\,\, $\nu\,(\,x_{n\,+\,p}\,-\,x_n
\,,\, t\,)\;<\;r$ \,\,for all $n\;\geq \;n_0$. Hence $\{\,x_{\,n}\}$ is a Cauchy sequence in $\left( {\,V\;,\;A^{\,\psi } \,} \right)$.

\begin{Definition}
A generalized IF$\psi$NLS
$\left( {\,V\;,\;A^{\,\psi } \,} \right)$ is said to be complete if every cauchy sequence in $\left( {\,V\;,\;A^{\,\psi } \,} \right)$ is convergent.
\end{Definition}

\begin{Theorem}
Let \,$\left( {\,V\;,\;A^{\,\psi } \,} \right)$\, be a generalized IF$\psi$NLS. A sufficient condition for the generalized IF$\psi$NLS \,$\left( {\,V\;,\;A^{\,\psi } \,} \right)$\, to be complete is that every Cauchy
sequence in \,$\left( {\,V\;,\;A^{\,\psi } \,} \right)$\, has a convergent subsequence.
\end{Theorem}

{\bf Proof.}$\;\;$ Let \,$\left\{ {\;x_{\,n}\;} \right\}_{\,n} $\, be a Cauchy
sequence in
$\left( {\,V\;,\;A^{\,\psi } \,} \right)$ and \,$\left\{ {\;x_{\,n_{\,k}}\;} \right\}_{\,k} $\ be a subsequence of \,$\left\{ {\;x_{\,n}\;} \right\}_{\,n} $\, that converges to $
 \,
 x\,\in  V $ and \,$s\,,\;t\,,\;s\;\circ\;t\, \;>\;
0 .$ \\Since \,$\left\{ {\;x_{\,n}\;} \right\}_{\,n} $\, is a cauchy sequence in $\left( {\,V\;,\;A^{\,\psi } \,} \right)$ , We have for $r \,>\, 0,\;r\,\in\,(\,0\,,\,1\,)$, there exists an integer $n_0\,\in\,N$ such that  \\ $\mu\,\left( {\,x_{\,n} \; - \;x_{\,k} \;,\; s\,}
\right)\; > \;1\;-\;r$ \,and\, $\nu\,\left( {\,x_{\,n} \; - \;x_{\,k} \;,\; s\,}
\right)\; < \;r$ \hspace{1.0cm}$\forall\;n \,,\,k\;
\geq\,n_0$  \\ Again , Since \,$\left\{ {\;x_{\,n}\;} \right\}_{\,k} $\,converges to x ,We have \\ $\mu\,\left( {\,x_{\,n_k} \; - \;x \;,\; t\,}
\right)\; > \;1 \;-\; r$ \,and\, $\nu\,\left( {\,x_{\,n_k} \; - \;x \;,\; t\,}
\right)\; < \;r$ \hspace{1.0cm}$\forall\;n_{\,k}\;
\geq\,n_0$ \\ Now,
\\ $\begin{array}{l}\mu\,\left( {\,x_{\,n} \; - \;x \;,\;s\;\circ\;t\,} \right)
 \;\; = \;\;\mu\,\left( {\,x_{\,n} \; - \,x_{\,n_{\,k}} \; + \,x_{\,n_{\,k}} \; - \;x
  \;,\;s\;\circ\;t\,} \right) \\
  \;\ge \;\; \mu\,\left( {\,x_{\,n} \; - \;x_{\,n_{\,k}}\;,\;s\,} \right)\;\; *
 \;\; \;\mu\,\left( {\,x_{\,n_{\,k}}
 \; - \;x\;,\;t\,} \right) \\
 \; > \;(\,1\;-\;r\,)\; * \;(\,1\;-\;r\,)\;\; = \;\;(\,1\;-\;r\,) \hspace{1.0cm}\forall\;n\;
\geq\,n_0 \
  \end{array}$
  \\Again , we see that $\\
\begin{array}{l}\nu\,\left( {\,x_{\,n} \; - \;x \;,\;s\;\circ\;t\,} \right)
 \;\; = \;\;\nu\,\left( {\,x_{\,n} \; - \,x_{\,n_{\,k}} \; + \,x_{\,n_{\,k}} \; - \;x
  \;,\;s\;\circ\;t\,} \right) \\
 \;\leq \;\; \nu\,\left( {\,x_{\,n} \; - \;x_{\,n_{\,k}}\;,\;s\,} \right)\;\; \diamond\
 \;\; \;\nu\,\left( {\,x_{\,n_{\,k}}
 \; - \;x\;,\;t\,} \right) \\
 \; < \;\;r\; \diamond\ \;r\;\; = \;\;r  \hspace{1.0cm}\forall\;n\;
\geq\,n_0\
  \end{array}$
\\Thus \,$\left\{ {\;x_{\,n}\;} \right\}_{\,n} $\, converges to x in $\left( {\,V\;,\;A^{\,\psi } \,} \right)$ and hence $\left( {\,V\;,\;A^{\,\psi } \,} \right)$ is complete.

\begin{Definition}
Let \,$\left( {\,V\;,\;A^{\,\psi } \,} \right)$\, be a generalized IF$\psi$NLS. A subset \,$P$\, of \,$V$\,
is said to be \textbf{closed} if for any sequence
\,$\{\,x_{\,n}\,\}_{\,n}$\, in \,$P$\, converges to \,$x
\;\,\varepsilon\;\, P$, that is,  \[ \mathop {\lim }\limits_{n\;
\to \;\infty }\;\mu\,(\,x_{\,n} \;-\; x \;,\; t\,) \;\,=\;\, 1  \;
and \;\mathop {\lim }\limits_{n\; \to \;\infty }\;\nu\,(\,x_{\,n} \;-\;
x \;,\; t\,) \;\,=\;\, 0 \;\; \Longrightarrow \;\; x
\;\,\varepsilon\;\, P .\]
\end{Definition}

\begin{Definition}
Let \,$\left( {\,V\;,\;A^{\,\psi } \,} \right)$\, be a generalized IF$\psi$NLS. A subset \,$Q$\, of \,$V$\,
is said to be the \textbf{closure} of \,$P \;(\;\subset \;V\;)$\, if
for any \,$x \;\,\varepsilon\;\, Q$ , there exists a sequence
\,$\{\,x_{\,n}\,\}_{\,n}$\, in \,$P$\, such that $\mathop{\lim }\limits_{n\; \to \;\infty } \;\mu\,\left( {\,x_{\,n} \; - \;x \;,\; t\,}
\right)\;\; = \;1 \;and\; \mathop {\lim }\limits_{n\; \to \;\infty
} \;\nu\,\left( {\,x_{\,n} \; - \;x \;,\; t\,}
\right)\; = \;0 \;\;\forall\;t\,\in\mathbf{R^{+}}.$  We denote the set
\,$Q$\, by \,$\overline{P}$
\end{Definition}

\begin{Definition}
A subset \,$P$\, of a generalized IF$\psi$NLS is said to be \textbf{bounded} if and
only if there exist \,$t \;>\; 0$\, and \,$0 \;<\; r \;<\; 1$\, such
that $\;\mu\,(\,x,\,t\,)\;>\;1\,-\,r$\,\,and\,\, $\nu\,(\,x,\,t\,)\;<\;r\;\;\;\forall\;x\,
\in\,P\;\,$.
\end{Definition}

\begin{Definition}
Let \,$\left( {\,V\;,\;A^{\,\psi } \,} \right)$\, be a generalized IF$\psi$NLS. A subset \,$P$\, of of
\,$V$\, is said to be \textbf{compact} if any sequence
\,$\{\,x_{\,n}\,\}_{\,n}$\, in \,$P$\, has a subsequence converging
to an element of \,$P$ .
\end{Definition}

Let \,$\left( {\,V\;,\;A^{\,\psi } \,} \right)$\, be a generalized IF$\psi$NLS. We further assume that \\\\ $(\,xii\,)$ \hspace{0.8cm} $\left. {{}_{a\;\; * \;\;a\;\; =
\;\;a}^{a\;\; \diamond \;\;a\;\; = \;\;a} \;\;}
\right\}\;\;\;\forall \;\;a\;\; \varepsilon \;\;[\,0\;\,,\;\,1\,]$ \\\\
$(\,xiii\,)$ \hspace{0.8cm} $\mu\,(\,x \;,\; t\,) \;>\; 0 \;\; \forall
\;\; t \;>\; 0 \;\; \Longrightarrow \;\; x \;=\; \underline{0}$ \\\\
$(\,xiv\,)$ \hspace{0.8cm} $\nu\,(\,x \;,\; t\,) \;<\; 1 \;\; \forall
\;\; t \;>\; 0 \;\; \Longrightarrow \;\; x \;=\; \underline{0}$ \\

\begin{Theorem}
Let \,$\left( {\,V\;,\;A^{\,\psi } \,} \right)$\, be a generalized IF$\psi$NLS satisfying the condition
\,$(\,Xii\,)$ . Every Cauchy sequence in \,$\left( {\,V\;,\;A^{\,\psi } \,} \right)$\, is
bounded .
\end{Theorem}

{\bf Proof.}$\;\;$ Let us consider a fixed \,$r_{\,0}$\, with \,$0 \;<\; r_{\,0} \;<\; 1$
\, and \,$\{\,x_{\,n}\,\}_{\,n}$\, be a Cauchy sequence in a generalized IF$\psi$NLS
\,$\left( {\,V\;,\;A^{\,\psi } \,} \right)$.  Then $\exists \;\; n_{\,0}\;\in\;\mathbf{N}$ such that \\\[
\left. {{}_{\;\nu\,\left( {\,x_{\,n\; + \;p} \; - \;x_{\,n} \;,\;t\,} \right)\;\;
< \;\;r_{\,0}}^{ \;\mu\,\left(
{\,x_{\,n\; + \;p} \; - \;x_{\,n} \;,\;t\,} \right)\;\; > \;1 \;-\; r_{\,0}}
\;\;} \right\}\;\;\forall \;\;t\;\, > \;\,0\;\,,\;\,p\;\; =
\;\;1\;\,,\;\,2\;\,,\;\; \cdots , \;\forall \;\; n \;\, > \;\, n_{\,0}\;.
\]\\
Now we see that  \[ \mu\,(\,x_{\,n\; + \;p} \;-\; x_{\,n} \;,\; t\,) \;\,>\;\, 1 \;-\; r_{\,0} \;\;\forall \;\;t\;\, > \;\,0\;\,,\;\,p\;\; = \;\;1  \;,\; 2
\;,\;\; \cdots \,,\; \;\forall \;\; n \;\, > \;\, n_{\,0}\]$\begin{array}{l}
 \Rightarrow \;\;\mu\,\left(-({\,x_{\,n} \; - \;x_{\,n\; + \;p}\;)\;,\;t\,}\, \right)\;\;\; > \;\;1 \;-\; r_{\,0} \\
  \Rightarrow \;\;\mu\,\left( {\,x_{\,n} \; - \;x_{\,n\; + \;p}\;,\;\frac{t}{\psi(\,\;-1\,)}\,} \right)\;\; > \;\;1 \;-\; r_{\,0}
  \\
  \Rightarrow \;\;\mu\,\left( {\,x_{\,n} \; - \;x_{\,n\; + \;p}\;,\frac{t}{\psi(\,\-1\,)}\,} \right)\hspace{0.7cm} > \;\;1 \;-\; r_{\,0}
  \\
  \Rightarrow \;\;\mu\,\left({\,x_{\,n} \; - \;x_{\,n\; + \;p}\;,\;t\,} \right)\hspace{1.3cm} > \;\;1 \;-\; r_{\,0} \\
 \end{array}$ \\\\
 $\;\; \Longrightarrow \;\; For \; t' \;>\; 0 \;\;
\exists \;\; n_{\,0} \;=\; n_{\,0}(\,t'\,)$ \\
such that $\mu\,(\,x_{\,n} \;-\; x_{\,n\; + \;p} \;,\; t'\,) \;\,>\;\, 1 \;-\; r_{\,0}
\;\; \forall \;\; n \;\, > \;\, n_{\,0} \;,\; p\;\; =
\;\;1\;\,,\;\,2\;\,,\;\; \cdots$ \\
Since \,$\mathop {\lim
}\limits_{t\; \to \;\infty }\;\mu\,(\,x \;,\; t\,) \;\,=\;\, 1$, we
have for each \,$x_{\,i}$ , \,$\exists \;\;t_{\,i} \;>\; 0$\, such
that \[\mu\,(\,x_{\,i} \;,\; t\,) \;>\; 1 \;-\; r_{\,0} \hspace{0.5cm}
\forall \;\; t \;\ge\; t_{\,i} \;,\; i \;=\; 1 \;,\; 2 \;,\; \cdots\]
Let \,$t_{\,0} \;=\; t' \;\circ\; \max\,\{\,t_{\,1} \;,\; t_{\,2}
\;,\; \cdots \;,\; t_{\,n_{\,0}} \,\}$ . Then , \\ $\mu\,(\,x_{\,n}
\;,\; t_{\,0}\,) \;\,\ge\;\, \mu\,(\,x_{\,n} \;,\; t' \;\circ\;
t_{\,n_{\,0}}\,) \\ {\hspace{2.2cm}} = \;\; \mu\,(\,x_{\,n} \;-\;
x_{\,n_{\,0}} \;+\; x_{\,n_{\,0}} \;,\; t' \;\circ\; t_{\,n_{\,0}}\,)
\\ {\hspace{2.2cm}} \ge \; \mu\,(\,x_{\,n} \;-\;
x_{\,n_{\,0}} \;,\; t'\,) \;\ast\; \mu\,(\,x_{\,n_{\,0}} \;,\;
t_{\,n_{\,0}}\,) \\{\hspace{2.2cm}} > \;\; (\,1 \;-\; r_{\,0}\,) \;\ast\; (\,1 \;-\; r_{\,0}\,)
\;\,=\;\, 1 \;-\; r_{\,0} \hspace{0.5cm} \forall \;\; n \;>\; n_{\,0}$ \\
Thus , we have \\ ${\hspace{2.5cm}} \mu\,(\,x_{\,n} \;,\; t_{\,0}\,)
\;\,>\;\, 1 \;-\; r_{\,0} \hspace{0.5cm} \forall \;\; n \;>\; n_{\,0}$ \\
Also , ${\hspace{0.2cm}} \mu\,(\,x_{\,n} \;,\; t_{\,0}\,) \;\,\ge\;\,
\mu\,(\,x_{\,n} \;,\; t_{\,n}\,) \;\,>\;\, 1 \;-\; r_{\,0} \hspace{0.5cm}\
\forall \;\;n \;=\; 1 \;,\; 2 \;,\; \cdots \;,\; n_{\,0}$ \\ So, we
have \\ ${\hspace{2.5cm}} \mu\,(\,x_{\,n} \;,\; t_{\,0}\,) \;\,>\;\,
1 \;-\; r_{\,0} \hspace{0.5cm} \forall \;\; n \;=\; 1 \;,\; 2 \;,\; \cdots
\hspace{0.8cm} \cdots \hspace{0.5cm} (\,1\,)$
\\ Again, we see that \\
$\nu\,(\,x_{\,n\; + \;p} \;-\; x_{\,n} \;,\; t\,) \;\,<\;\,r_{\,0} \;\;\forall
\;\;t\;\, > \;\,0\;\,,\;\,p\;\; = \;\;1 \;,\; 2 \;,\;\; \cdots \;\forall \;\; n \;>\; n_{\,0}$  $\begin{array}{l}
 \Rightarrow \; \nu\,\left(-({\,x_{\,n} \; - \;x_{\,n\; + \;p}\;)\;,\;t\,} \,\right)\hspace{0.47cm} < \;\;r_{\,0} \\
  \Rightarrow \;\nu\,\left( {\,x_{\,n} \; - \;x_{\,n\; + \;p}\;,\;\frac{t}{\psi(\,\;-1\,)}\,} \,\right)\;\, < \;\;r_{\,0} \\
\Rightarrow \;\nu\,\left( {\,x_{\,n} \; - \;x_{\,n\; + \;p}\;,\;\frac{t}{\psi(\,\-1\,)}\,} \right)\hspace{0.54cm} < \;\;r_{\,0}
\\ \Rightarrow \;\nu\,\left({\,x_{\,n} \; - \;x_{\,n\; + \;p}\;,\;t\,}\, \right)\hspace{1.15cm} < \;\;r_{\,0} \\
 \end{array}$
 \\ \\
 $\Longrightarrow \;\; For \; t' \;>\; 0 \;\; \exists \;\; n'_{\,0}
\;=\; n'_{\,0}(\,t'\,)$ \, such that  \\$\;\nu\,(\,x_{\,n} \;-\;
x_{\,n\; + \;p} \;,\; t'\,) \;\,<\;\, r_{\,0} \;\;
\forall \;\; n \;\, > \;\, n'_{\,0} \;,\; p\;\; =
\;\;1\;\,,\;\,2\;\,,\;\; \cdots$
\\ Since \,$\mathop {\lim
}\limits_{t\; \to \;\infty }\;\nu\,(\,x \;,\; t\,) \;\,=\;\, 0$, we
have for each \,$x_{\,i}$ ,  $\exists \;t'_{\,i} \;>\; 0$ such
that \[\nu\,(\,x_{\,i} \;,\; t\,) \;<\; r_{\,0}
\hspace{0.5cm} \forall \;\; t \;\ge\; t'_{\,i} \;,\; i \;=\; 1 \;,\; 2
\;,\; \cdots\]  Let \,$t'_{\,0} \;=\; t' \circ\ \max\,\{\,t'_{\,1}
\;,\; t'_{\,2} \;,\; \cdots \;,\; t'_{\,n_{\,0}} \,\}$ .Then , \\
$\nu\,(\,x_{\,n} \;,\; t'_{\,0}\,) \hspace{0.37cm} \le\;\, \nu\,(\,x_{\,n} \;,\; t'
\;\circ\; t'_{\,n_{\,0}}\,) \\ {\hspace{2.4cm}} = \;\; \nu\,(\,x_{\,n}
\;-\; x_{\,n'_{\,0}} \;+\; x_{\,n'_{\,0}} \;,\; t' \;\circ\;
t'_{\,n_{\,0}}\,)\\ {\hspace{2.4cm}} \le \;\; \nu\,(\,x_{\,n} \;-\;
x_{\,n'_{\,0}} \;,\; t'\,) \;\diamond\; \nu\,(\,x_{\,n'_{\,0}} \;,\;
t'_{\,n_{\,0}}\,) \\{\hspace{2.5cm}} < \;\; r_{\,0}
\;\diamond\; r_{\,0} \;\,=\;\, r_{\,0}
\hspace{0.5cm} \forall \; n \;>\; n'_{\,0}$ \\ Thus , we have \\
${\hspace{2.5cm}} \nu\,(\,x_{\,n} \;,\; t'_{\,0}\,) \;\,<\;\, r_{\,0} \hspace{0.5cm} \forall \;\; n \;>\; n'_{\,0}$ \\
Also , $\nu\,(\,x_{\,n} \;,\; t'_{\,0}\,) \;\,\le\;\,
\nu\,(\,x_{\,n} \;,\; t'_{\,n}\,) \;\,<\;\, r_{\,0}
\hspace{0.5cm}\forall \;\;n \;=\; 1
\;,\; 2 \;,\; \cdots \;,\; n'_{\,0}$ \\ So, we have \\
${\hspace{2.5cm}} \nu\,(\,x_{\,n} \;,\; t'_{\,0}\,) \;\,<\;\, r_{\,0} \hspace{0.5cm} \forall \;\; n \;=\; 1 \;,\; 2 \;,\;
\cdots \hspace{0.8cm} \cdots \hspace{0.5cm} (\,2\,)$ \\ Let
$t''_{\,0} \;=\; \max\,\{\,t_{\,0} \;,\; t'_{\,0}\,\}$ . Hence from
\,$(\,1\,)$\, and \,$(\,2\,)$\, we see that \\
\[
\left. {{}_{\nu\,\left( {\,x_{\,n} \;,\;t''_{\,0} \,} \right)\;\, <
\;\;r_{\,0} }^{\mu\,\left( {\,x_{\,n}
\;,\;t''_{\,0} \,} \right)\;\;\, > \;\;\left( {\,1\; - \;r_{\,0} \,} \right) } \;\;}
\right\}\;\;\forall \;\;n\;\; = \;\;1\;\,,\;\,2\;\,,\;\; \cdots
\] \\ This implies that \,$\{\,x_{\,n}\,\}_{\,n}$\,is
bounded in \,$\left( {\,V\;,\;A^{\,\psi } \,} \right)$

\begin{Theorem}
In a finite dimensional generalized IF$\psi$NLS \,$\left( {\,V\;,\;A^{\,\psi } \,} \right)$\, satisfying the
conditions \,$(\,Xii\,)$, \,$(\,Xiii\,)$\, and \,$(\,Xiv\,)$ , a
subset \,$P$\, of \,$V$\, is compact if and only if \,$P$\, is
closed and bounded in \,$(\,V \;,\; A\,)$.
\end{Theorem}
{\bf Proof.}
$\Longrightarrow \;\; part \; :\;$ Proof of this part directly
follows from the proof of the theorem 2.5 \cite{Bag1}. \\
$\Longleftarrow \; part \; : \;$ In this part, we suppose that
\,$P$\, is closed and bounded in the finite dimensional generalized IF$\psi$NLS
\,$\left( {\,V\;,\;A^{\,\psi } \,} \right)$. To show \,$P$\, is compact, consider
\,$\{\,x_{\,n}\,\}_{\,n}$, an arbitrary sequence in \,$P$. Since
\,$V$\, is finite dimensional, let \,$\dim\,V \;=\; n$\, and
\,$\{\,e_{\,1} \;,\; e_{\,2} \;,\; \cdots \;,\; e_{\,n} \,\}$\, be a
basis of \,$V$. So, for each \,$x_{\,k}$, \,$\exists \;\;
\beta_{\,1}^{\,k} \;,\; \beta_{\,2}^{\,k} \;,\; \cdots \;,\;
\beta_{\,n}^{\,k} \;\; \varepsilon \;\; F$\, such that
 \[x_{\,k} \;\,=\;\, \beta_{\,1}^{\,k}\;e_{\,1} \;+\; \beta_{\,2}^{\,k}\;
e_{\,2} \;+\; \cdots \;+\; \beta_{\,n}^{\,k}\;e_{\,n} \;,\; k \;=\;
1 \;,\; 2 \;,\; \cdots\]  
Following the calculation of the theorem 11 \cite{Samanta}, we can write\\ \,$x_{\,k_{\,l}} \;=\; \beta_{\,1}^{\,k_{\,l}}\;e_{\,1} \;+\;
\beta_{\,2}^{\,k_{\,l}}\;e_{\,2} \;+\; \cdots \;+\;
\beta_{\,n}^{\,k_{\,l}}\;e_{\,n}$\, and \,$\beta_{\,1} \;=\; \mathop
{\lim }\limits_{n\; \to \;\infty }\;\beta_{\,1}^{\,k_{\,l}} \;,\;
\beta_{\,2} \;=\; \mathop {\lim }\limits_{n\; \to \;\infty
}\;\beta_{\,2}^{\,k_{\,l}} \;,\; \cdots \;,\; \beta_{\,n} \;=\;
\mathop {\lim }\limits_{n\; \to \;\infty
}\;\beta_{\,n}^{\,k_{\,l}}$\, and \,$x \;=\; \beta_{\,1}\,e_{\,1}
\;+\; \beta_{\,2}\,e_{\,2} \;+\; \cdots \;+\; \beta_{\,n}\,e_{\,n}$.
\\ Now suppose that for all \, $t \;>\; 0$, there exist $t_{\,1} \,,\,
t_{\,2} \,,\; \cdots \; , t_{\,k} \;>\; 0 $ \,such that\,
$ (\,t_{\,1}\,\circ
\, t_{\,2}\,\circ \; \cdots \; \circ \,t_{\,k} \,)\; >\;0$. Then we have \\
$\mu\,(\,x_{\,k_{\,l}} \;-\; x\, \;,\; t_{\,1}\,\circ
\, t_{\,2}\,\circ \; \cdots \; \circ \,t_{\,k}\,) 
\\{\hspace{3.1cm}}=\;\, \mu\,(\,
\sum\limits_{i\; = \;1}^n {\beta _{\,i}^{\,k_{\,l}} \,e_{\,i} }
\;-\; \sum\limits_{i\; = \;1}^n {\beta _{\,i} \,e_{\,i} } \;,\; t_{\,1}\,\circ
\, t_{\,2}\,\circ \; \cdots \; \circ \,t_{\,k} \,) \\
{\hspace{3.1cm}}=\;\; \mu\,(\, \sum\limits_{i\; = \;1}^n \,(\,{\beta
_{\,i}^{\,k_{\,l}} \;-\; \beta_{\,i}\,) \,e_{\,i} } \;,\; t_{\,1}\,\circ
\, t_{\,2}\,\circ \; \cdots \; \circ \,t_{\,k}\,) \\
{\hspace{3.1cm}} \ge \;\; \mu\,(\,(\,\beta _{\,1}^{\,k_{\,l}} \;-\;
\beta_{\,1} \,)\;e_{\,1} \;,\; t_{\,1} \;) \;\ast\; \cdots
\;\ast\; \mu\,(\,(\,\beta _{\,n}^{\,k_{\,l}} \;-\; \beta_{\,n}
\,)\;e_{\,n} \;,\; t_{\,n} \;) \\ {\hspace{3.1cm}} = \;\;
\mu\,\left(\,e_{\,1} \;,\; \frac{t_{\,1}}{\psi\,\left(\,\beta _{\,1}^{\,k_{\,l}} \;-\;
\beta_{\,1} \,\right)} \,\right) \;\ast\; \cdots \;\ast\; \mu\,\left(\,e_{\,n} \;,\;
\frac{t_{\,n}}{\psi\,\left(\,\beta _{\,n}^{\,k_{\,l}} \;-\; \beta_{\,n} \,\right)} \,\right)$
\\\\ Since \, $\mathop {\lim }\limits_{l\; \to \;\infty }
\;\frac{t_{\,i}}{{\psi\,\left( {\,\beta _i^{\,k_{\,l}} \; - \;\beta _{\,i}
\,} \right)}}\;\, = \;\,\infty $, we see that \\${\hspace{4.5cm}}\mathop {\lim
}\limits_{l\; \to \;\infty }\;\mu\,\left(\,e_{\,i} \;,\;
\frac{t_{\,i}}{{\psi\,\left( {\,\beta _i^{\,k_{\,l}} \; - \;\beta _{\,i}
\,}\right)}} \,\right)
\;=\; 1$ \\\\ $\Longrightarrow\;\; \mathop {\lim }\limits_{l\; \to
\;\infty }\;\mu\,(\,x_{\,k_{\,l}} \;-\; x \;,\; t\,) \;\,\ge\;\, 1
\;\ast\; \cdots \;\ast\; 1 \;=\; 1 \hspace{0.5cm} \forall \;\; t \;>\; 0$ \\
$\Longrightarrow\;\; \mathop {\lim }\limits_{l\; \to \;\infty
}\;\mu\,(\,x_{\,k_{\,l}} \;-\; x \;,\; t\,) \;\,=\;\, 1 \hspace{0.5cm}
\forall \;\; t \;>\; 0 \hspace{1.5cm} \cdots \hspace{1.5cm} (\,4\,)$
\\ Again, we have \\
$\nu\,(\,x_{\,k_{\,l}} \;-\; x\, \;,\; t_{\,1}\,\circ
\, t_{\,2}\,\circ \; \cdots \; \circ \,t_{\,k}\,) \\ 
{\hspace{3.1cm}}=\;\, \nu\,(\,
\sum\limits_{i\; = \;1}^n {\beta _{\,i}^{\,k_{\,l}} \,e_{\,i} }
\;-\; \sum\limits_{i\; = \;1}^n {\beta _{\,i} \,e_{\,i} } \;,\; t_{\,1}\,\circ
\, t_{\,2}\,\circ \; \cdots \; \circ \,t_{\,k} \,) \\
{\hspace{3.1cm}}=\;\; \nu\,(\, \sum\limits_{i\; = \;1}^n \,(\,{\beta
_{\,i}^{\,k_{\,l}} \;-\; \beta_{\,i}\,) \,e_{\,i} } \;,\; t_{\,1}\,\circ
\, t_{\,2}\,\circ \; \cdots \; \circ \,t_{\,k} \,) \\
{\hspace{3.1cm}} \le \;\; \nu\,(\,(\,\beta _{\,1}^{\,k_{\,l}} \;-\;
\beta_{\,1} \,)\;e_{\,1} \;,\; t_{\,1} \;) \;\diamond\; \cdots
\;\diamond\; \nu\,(\,(\,\beta _{\,n}^{\,k_{\,l}} \;-\; \beta_{\,n}
\,)\;e_{\,n} \;,\; t_{\,n} \;) \\ {\hspace{3.1cm}} = \;\;
\nu\,\left(\,e_{\,1} \;,\; \frac{t_{\,1}}{\psi\,\left(\,\beta _{\,1}^{\,k_{\,l}} \;-\;
\beta_{\,1} \,\right)} \,\right) \;\diamond\; \cdots \;\diamond\; \nu\,\left(\,e_{\,n}
\;,\; \frac{t_{\,n}}{\psi\,\left(\,\beta _{\,n}^{\,k_{\,l}} \;-\; \beta_{\,n} \,\right)}
\,\right)$ \\\\ Since \, $\mathop {\lim }\limits_{l\; \to \;\infty }
\;\frac{t_{\,i}}{{\psi\,\left( {\,\beta _i^{\,k_{\,l}} \; - \;\beta _{\,i}
\,} \right)}}\;\, = \;\,\infty $, we see that \\${\hspace{4.5cm}}\mathop {\lim
}\limits_{l\; \to \;\infty }\;\nu\,\left(\,e_{\,i} \;,\;
\frac{t_{\,i}}{{\psi\,\left( {\,\beta _i^{\,k_{\,l}} \; - \;\beta _{\,i}
\,}\right)}} \,\right)
\;=\; 0$ \\\\ $\Longrightarrow\;\; \mathop {\lim }\limits_{l\; \to
\;\infty }\;\nu\,(\,x_{\,k_{\,l}} \;-\; x \;,\; t\,) \;\,\le\;\, 0
\;\diamond\; \cdots \;\diamond\; 0 \;=\; 0 \hspace{0.5cm} \forall
\;\; t \;>\; 0$ \\ $\Longrightarrow\;\; \mathop {\lim }\limits_{l\;
\to \;\infty }\;\nu\,(\,x_{\,k_{\,l}} \;-\; x \;,\; t\,) \;\,=\;\, 0
\hspace{0.5cm} \forall \;\; t \;>\; 0 \hspace{1.5cm} \cdots
\hspace{1.5cm} (\,5\,)$ \\ Thus, from \,$(\,4\,)$\, and
\,$(\,5\,)$\, we see that \\${\hspace{2.5cm}} \mathop {\lim
}\limits_{l\; \to \;\infty }\,x_{\,k_{\,l}} \;\,=\;\, x \;\;
\Longrightarrow \;\; x \;\,\varepsilon \;\, A$ \, $[$\, Since
\,$A$\, is closed \,$]$. \\${\hspace{2.5cm}} \Longrightarrow \;\; A
\;\; is \;\; compact$.


\section{ Generalized Intuitionistic Fuzzy \\ $\psi-\alpha$-Normed Linear space}


\begin{Theorem}
Define \,$\left\| {\;x\;} \right\|_{\,\alpha }^{\,1}  \;\,=\;\,
\wedge\;\{\,t \;\,:\;\, \mu\,(\,x \;,\; t\,) \;\,\ge\;\, \alpha
\,\}$\, and \,$\left\| {\;x\;} \right\|_{\,\alpha }^{\,2}  \;\,=\;\,
\vee\;\{\,t \;\,:\;\, \nu\,(\,x \;,\; t\,) \;\,\le\;\, \alpha \,\}$\,
, $\alpha \;\,\varepsilon\;\, (\,0 \;,\; 1\,)$ . Then both
\,$\{\;\left\| {\;x\;} \right\|_{\,\alpha }^{\,1} \;\;:\;\; \alpha
\;\,\varepsilon\;\, (\,0 \;,\; 1\,) \;\}$\, and \,$\{\;\left\|
{\;x\;} \right\|_{\,\alpha }^{\,2} \;\;:\;\; \alpha
\;\,\varepsilon\;\, (\,0 \;,\; 1\,) \;\}$\, are ascending family of
norms on \,$V$ . These norms are said to be \,$\alpha$ - norm on \,$V$\,
corresponding to the generalized IF$\psi$NLS \,$A$\, on \,$V$ .
\end{Theorem}

{\bf Proof.}$\;\;$
Let \,$\alpha \;\,\varepsilon\;\, (\,0 \;,\; 1\,)$ . To prove
\,$\left\| {\;x\;} \right\|_{\,\alpha }^{\,1}$\, is a norm on \,$V$
, we will prove the followings $:$\\\\ $(\,1\,)$ \hspace{1.5cm}
$\left\| {\;x\;} \right\|_{\,\alpha }^{\,1} \;\,\ge\;\, 0$
\hspace{0.8cm} $\forall \;\; x \;\; \varepsilon \;\; V$ ; \\\\
$(\,2\,)$ \hspace{1.5cm} $\left\| {\;x\;} \right\|_{\,\alpha }^{\,1}
\;\,=\;\, 0 \;\;\Longleftrightarrow\;\; x \;=\; \underline{0}$ ;
\\\\ $(\,3\,)$ \hspace{1.5cm} $\left\| {\;\alpha\;x\;} \right\|_{\,\alpha
}^{\,1} \;\,=\;\, $  $\psi (\,\alpha \,)\;\left\| {\;x\;} \right\|_{\,\alpha
}^{\,1}$ ;\\\\ $(\,4\,)$ \hspace{1.5cm} $\left\| {\;x \;+\; y\;}
\right\|_{\,\alpha }^{\,1} \;\,\le\;\, \left\| {\;x\;}
\right\|_{\,\alpha }^{\,1} \;+\; \left\| {\;y\;} \right\|_{\,\alpha
}^{\,1}$ .\\\\ The proof of $(\,1\,)$ and $(\,2\,)$
directly follows from the proof of the theorem 2.1 \cite{Bag1} . So,
we now prove $(\,3\,)$ and $(\,4\,)$ .\\If $\alpha\;\;=\;\;0$ \,and\, $\psi (\,\alpha \,)\;\;=\;\;|\,\alpha\,|$ then\\ $\left\| {\;\alpha\;x\;} \right\|_{\,\alpha
}^{\,1} \;\,=\;\, $  $\left\| {\; \underline{0}\;} \right\|_{\,\alpha
}^{\,1}\;\,=\;\, 0 \;\,=\;\, 0 \;\left\| {\;x\;} \right\|_{\,\alpha
}^{\,1}\;\,=\;\,|\,\alpha \,|\left\| {\;x\;} \right\|_{\,\alpha
}^{\,1}\;\,=\;\,\psi (\,\alpha \,)\;\left\| {\;x\;} \right\|_{\,\alpha
}^{\,1}$\\If $\alpha\;\;\neq\;\;0$ then\\ $\left\| {\;\alpha\;x\;} \right\|_{\,\alpha
}^{\,1}\\ \;\,=\;\,\wedge\;\{\,s \;\,:\;\, \mu\,(\;\alpha\;x\;,\; s\,) \;\,\ge\;\, \alpha \,\}\\ \;\,=\;\,\wedge\;\{\,s \;\,:\;\, \mu\,\left( { \;x\;,\frac{s}{\psi(\,\alpha\,)}\,} \right) \;\,\ge\;\, \alpha \,\}\\\;\; = \;\,\wedge\;\{\psi (\,\alpha \,)s \;\,:\;\, \mu\,(\;\;x\;,\; s\,) \;\,\ge\;\, \alpha \,\}\\ \;\,=\;\,\wedge\;\psi (\,\alpha \,)
\;\{s \;\,:\;\, \mu\,(\;\;x\;,\; s\,) \;\,\ge\;\, \alpha \,\}$
\\Therfore $\left\| {\;\alpha\;x\;} \right\|_{\,\alpha
}^{\,1} \;\,=\;\,\psi (\,\alpha \,)\;\left\| {\;x\;} \right\|_{\,\alpha
}^{\,1}$
       \\\\ $\left\| {\;x\;} \right\|_{\,\alpha
}^{\,1} \;+\; \left\| {\;y\;} \right\|_{\,\alpha }^{\,1}\\ \,\;=\,\;
\wedge\;\{\,s \;\,:\;\, \mu\,(\,x \;,\; s\,) \;\,\ge\;\, \alpha \,\}
\;+\; \wedge\;\{\,t \;\,:\;\, \mu\,(\,y \;,\; t\,) \;\,\ge\;\, \alpha
\,\}\\ \;\,\ge\;\, \wedge\;\{\,s \;\circ\; t \;\,:\;\, \mu\,(\,x \;,\; s\,)
\;\,\ge\;\, \alpha  \;,\; \mu\,(\,y \;,\; t\,) \;\,\ge\;\, \alpha\,\}\\
\;\,=\;\, \wedge\;\{\,s \;\circ\; t\;\,:\;\, \mu\,(\,x \;,\; s\,)
\;\ast\; \mu\,(\,y \;,\; t\,) \;\,\ge\;\, \alpha \;\ast\; \alpha \,\}\\
\;\,\ge\;\, \wedge\;\{\,s \;\circ\; t \;\,:\;\, \mu\,(\,x \;+\; y \;,\; s
\;\circ\; t\,) \;\,\ge\;\, \alpha \,\}\\ \;\,=\;\, \left\| {\;x \;+\; y\;}
\right\|_{\,\alpha }^{\,1}$, which proves \,$(\,4\,)$ .\\ Let \,$0
\;\,<\;\, \alpha_{\,1} \;\,<\;\, \alpha_{\,2} \;\,<\;\, 1$ .\\$\left\| {\;x\;} \right\|_{\,\alpha_{\,1} }^{\,1} \;\,=\;\,
\wedge\;\{\,t \;\,:\;\, \mu\,(\,x \;,\; t\,) \;\,\ge\;\, \alpha_{\,1}
\,\}$\,\\ and \,$\left\| {\;x\;} \right\|_{\,\alpha_{\,2} }^{\,1}
\;\,=\;\, \wedge\;\{\,t \;\,:\;\, \mu\,(\,x \;,\; t\,) \;\,\ge\;\,
\alpha_{\,2} \,\}$ .\\ Since \,$\alpha_{\,1} \;\,<\;\, \alpha_{\,2}$
$\{\,t \;\,:\;\, \mu\,(\,x \;,\; t\,) \;\,\ge\;\, \alpha_{\,2} \,\}
\;\;\subset\;\; \{\,t \;\,:\;\, \mu\,(\,x \;,\; t\,) \;\,\ge\;\,
\alpha_{\,1} \,\}\\ \;\;\Longrightarrow\;\; \wedge\{\,t \;\,:\;\,
\mu\,(\,x \;,\; t\,) \;\,\ge\;\, \alpha_{\,2} \,\} \;\;\ge\;\;
\wedge\{\,t \;\,:\;\, \mu\,(\,x \;,\; t\,) \;\,\ge\;\, \alpha_{\,1}
\,\}\\ \;\;\Longrightarrow\;\; \left\| {\;x\;}
\right\|_{\,\alpha_{\,2} }^{\,1} \;\,\le\;\, \left\| {\;x\;}
\right\|_{\,\alpha_{\,1} }^{\,1}$ .Thus, we see that \,$\{\;\left\|
{\;x\;} \right\|_{\,\alpha }^{\,1} \;\;:\;\; \alpha
\;\,\varepsilon\;\, (\,0 \;,\; 1\,) \;\}$\, is an ascending family
of norms on \,$V$ . \\ Now we shall prove that \,$\{\;\left\|
{\;x\;} \right\|_{\,\alpha }^{\,2} \;\;:\;\; \alpha
\;\,\varepsilon\;\, (\,0 \;,\; 1\,) \;\}$\, is also an ascending
family of norms on \,$V$.\\ Let \,$\alpha \;\,\varepsilon\;\, (\,0
\;,\; 1\,)$\, and \,$x \;,\; y \;\,\varepsilon\;\, V$ .\\It is
obvious that \,$\left\| {\;x\;} \right\|_{\,\alpha }^{\,2}
\;\,\le\;\, 0$ .\\Let \,$\left\| {\;x\;} \right\|_{\,\alpha }^{\,2}
\;\,=\;\, 0$ . Now,\\ \,$\left\| {\;x\;} \right\|_{\,\alpha }^{\,2}
\;\,=\;\, 0  \;\;\Longrightarrow\;\; \wedge\{\,t \;\,:\;\,\nu\,(\,x
\;,\; t\,) \;\,\le\;\, \left( {\,1\; - \alpha
 \,} \right)\} \;\,=\;\, 0
\\ \;\;\Longrightarrow\;\; \nu\,(\,x \;,\; t\,) \;\,>\;\, \alpha \;>\; 0
\;\; \forall \;\; t \;>\; 0 \;\;\Longrightarrow\;\; x \;=\;
\underline{0}$ .\\ Conversely, we assume that \,$x \;=\; \underline{0}
\;\;\Longrightarrow\;\;\nu\,(\,x \;,\; t\,) \;\,=\;\, 0 \;\; \forall
\;\; t \;>\; 0 \\ \;\;\Longrightarrow\;\; \wedge\{\,t \;\,:\;\,\nu\,(\,x
\;,\; t\,) \;\,\le\;\, \left( {\,1\; - \alpha
 \,} \right)\} \;\,=\;\, 0
\;\;\Longrightarrow\;\; \left\| {\;x\;} \right\|_{\,\alpha }^{\,2}
\;\,=\;\, 0$ .\\If $\alpha\;\;=\;\;0$ \, and \, $\psi (\,\alpha \,)\;\;=\;\;|\,\alpha\,|$ then \\ $\left\| {\;\alpha\;x\;} \right\|_{\,\alpha
}^{\,1} \;\,=\;\, $  $\left\| {\; \underline{0}\;} \right\|_{\,\alpha
}^{\,1}\;\,=\;\, 0 \;\,=\;\, 0 \;\left\| {\;x\;} \right\|_{\,\alpha
}^{\,1}\;\,=\;\,|\,\alpha \,|\left\| {\;x\;} \right\|_{\,\alpha
}^{\,1}\;\,=\;\,\psi (\,\alpha \,)\;\left\| {\;x\;} \right\|_{\,\alpha
}^{\,1}$\\If $\alpha\;\;\neq\;\;0$ \, then \\ $\left\| {\;\alpha\;x\;} \right\|_{\,\alpha
}^{\,1}\\ \;\,=\;\,\wedge\;\{\,s \;\,:\;\, \nu\,(\;\alpha\;x\;,\; s\,) \;\,\le\;\, \alpha
  \,\}\\ \;\,=\;\,\wedge\;\{\,s \;\,:\;\, \nu\,\left( { \;x\;,\frac{s}{\psi(\,\alpha\,)}\,} \right) \;\,\le\;\, \alpha
 \}\\\;\; = \;\,\wedge\;\{\psi (\,\alpha \,)s \;\,:\;\, \nu\,(\,x\;,\; s\,) \;\,\le\;\, \alpha \}\\ \;\,=\;\,\wedge\;\psi (\,\alpha \,)\;\{s \;\,:\;\, \nu\,(\,x\;,\; s\,) \;\,\le\;\, \alpha \}$
 \\Therfore $\left\| {\;\alpha\;x\;} \right\|_{\,\alpha
}^{\,1} \;\,=\;\,\psi (\,\alpha \,)\;\left\| {\;x\;} \right\|_{\,\alpha
}^{\,1}$. \\ $\left\| {\;x\;} \right\|_{\,\alpha
}^{\,2} \;+\; \left\| {\;y\;} \right\|_{\,\alpha }^{\,2}\\\;\;=\,\;
\wedge\;\{\,s \;\,:\;\, \nu\,(\,x \;,\; s\,) \;\,\le\;\, \alpha \}
\;+\; \wedge\;\{\,t \;\,:\;\, \nu\,(\,y \;,\; t\,) \;\,\le\;\, \alpha
\,\} \\\;\;\le\;\, \wedge\;\{\,s \;\circ\; t \;\,:\;\, \nu\,(\,x \;,\; s\,)
\;\,\le\;\, \alpha
  \,\,,\; \nu\,(\,y \;,\; t\,) \;\,\le\;\, \alpha \}
\\\;\;=\;\, \wedge\;\{\,s \;\circ\; t \;\,:\;\, \nu\,(\,x \;,\; s\,)
\;\diamond\; \nu\,(\,y \;,\; t\,) \;\,\le\;\, \alpha
  \,\diamond\; \alpha
 \,\} \\\;\,\le\;\, \wedge\;\{\,s \;\circ\; t \;\,:\;\, \nu\,(\,x \;+\;
y \;,\; s \;\circ\; t\,) \;\,\le\;\,\alpha
 \}  \\\;\,=\;\, \left\| {\;x
\;+\; y\;} \right\|_{\,\alpha }^{\,2}$ ,\\ That is \,$\left\| {\;x
\;+\; y\;} \right\|_{\,\alpha }^{\,2} \;\,\le\;\, \left\| {\;x\;}
\right\|_{\,\alpha }^{\,2} \;+\; \left\| {\;y\;} \right\|_{\,\alpha
}^{\,2}$\hspace{0.3cm} $\forall \;\;x \;,\; y \;\,\varepsilon\;\,
V$. \\ Let \,$0 \;\,<\;\, \alpha_{\,1} \;\,<\;\, \alpha_{\,2}
\;\,<\;\, 1$ .\\ Therefore, \,$\left\| {\;x\;}
\right\|_{\,\alpha_{\,1} }^{\,2} \;\,=\;\, \wedge\;\{\,t \;\,:\;\,
\nu\,(\,x \;,\; t\,) \;\,\le\;\, \alpha_{\,1}
  \,\}$ \,and\,\\ $\left\|
{\;x\;} \right\|_{\,\alpha_{\,2} }^{\,2} \;\,=\;\, \wedge\;\{\,t
\;\,:\;\, \nu\,(\,x \;,\; t\,) \;\,\le\;\,  \alpha_{\,2}
  \,\}$ . Since
\,$\alpha_{\,1} \;<\; \alpha_{\,2}$ , we have \\ $\{\,t \;\,:\;\,
\nu\,(\,x \;,\; t\,) \;\,\le\;\,  \alpha_{\,1}
  \,\} \;\;\subset\;\;
\{\,t \;\,:\;\, \nu\,(\,x \;,\; t\,) \;\,\le\;\, \alpha_{\,2}
  \,\}$ \\
$\Longrightarrow \;\;\; \wedge\{\,t \;\,:\;\,\nu\,(\,x \;,\; t\,)
\;\,\le\;\,  \alpha_{\,1}
  \,\} \;\;\le\;\; \wedge\{\,t \;\,:\;\,
\nu\,(\,x \;,\; t\,)\;\,\le\;\, \alpha_{\,2}
  \,\}$ \\
$\Longrightarrow \;\;\; \left\| {\;x\;} \right\|_{\,\alpha_{\,1}
}^{\,2} \;\,\le\;\, \left\| {\;x\;} \right\|_{\,\alpha_{\,2}
}^{\,2}$.\\  Thus we see that \,$\{\;\left\| {\;x\;} \right\|_{\,\alpha
}^{\,2} \;\;:\;\; \alpha \;\,\varepsilon\;\, (\,0 \;,\; 1\,) \;\}$\,
is an ascending family of norms on \,$V$ .

\begin{Lemma}\cite{Bag1}
\, Let  \,$\left( {\,V\;,\;A^{\,\psi } \,} \right)$\, be a generalized IF$\psi$NLS satisfying the
condition \,$(\,Xiii\,)$\, and \,$\{\,x_{\,1} \;,\; x_{\,2} \;,\;
\cdots \;,\; x_{\,n} \,\}$\, be a finite set of linearly independent
vectors of \,$V$ . Then for each \,$\alpha \;\,\varepsilon\;\, (\,0
\;,\; 1\,)$\, there exists a constant \,$C_{\,\alpha} \;>\; 0$\,
such that for any scalars \, $\alpha_{\,1} \;,\; \alpha_{\,2} \;,\;
\cdots \;,\; \alpha_{\,n}$ ,
\\ \[\|\,\alpha_{\,1}\,x_{\,1} \;+\; \alpha_{\,2}\,x_{\,2} \;+\;
\cdots \;+\; \alpha_{\,n}\,x_{\,n} \,\|_{\,\alpha}^{\,1} \;\,\ge\;\,
C_{\,\alpha}\;\sum\limits_{i\; = \;1}^n {\left| {\;\alpha _{\,i} \;}
\right|}\] \\ where \,$\|\,\cdot\,\|_{\,\alpha}^{\,1}$\, is defined
in the previous theorem.
\end{Lemma}

\begin{Theorem}
Every finite dimensional in a generalized IF$\psi$NLS satisfying the conditions
\,$(\,Xii\,)$\, and \,$(\,Xiii\,)$\, is complete .
\end{Theorem}

\medskip

{\bf Proof.}$\;\;$
Let \,$\left( {\,V\;,\;A^{\,\psi } \,} \right)$\, be a finite dimensional generalized IF$\psi$NLS satisfying
the conditions \,$(\,Xii\,)$\, and \,$(\,Xiii\,)$. Also, let
\,$\dim\,V \,=\, k$\, and \,$\{\,e_{\,1} \,,\, e_{\,2} \,,\; \cdots
\;,\, e_{\,k}\,\}$\, be a basis of \,$V$. Consider
\,$\{\,x_{\,n}\,\}_{\,n}$\, as an arbitrary Cauchy sequence in
\,$\left( {\,V\;,\;A^{\,\psi } \,} \right)$. \\ Let \,$x_{\,n} \;\,=\;\,
\beta_{\,1}^{\,(\,n\,)}\;e_{\,1} \;+\;
\beta_{\,2}^{\,(\,n\,)}\;e_{\,2} \;+\; \cdots \;+\;
\beta_{\,k}^{\,(\,n\,)}\;e_{\,k}$\, where \,$\beta_{\,1}^{\,(\,n\,)}
\;,\; \beta_{\,2}^{\,(\,n\,)} \;,\; \cdots \;,\;
\beta_{\,k}^{\,(\,n\,)}$\, are suitable scalars.
Now $\exists \; \beta_{\,1}\,,\, \beta_{\,2}\,,\cdots \; , \;\beta_{\,k} \,
 \in \;F$  ( By the calculation of the theorem 2.4 \cite{Bag1} )\\ Let \,
 $ x\;=\;\sum\limits_{i\; = \;1}^k
\,\beta_{\,i}\,e_{\,i}$, \,clearly \, $x\in\;V$.
Suppose that for all $t \;>\; 0$, there exist $t_{\,1} \,,\,
t_{\,2} \,,\; \cdots \; , t_{\,k} \;>\; 0 $ \,such that\,
$ (\,t_{\,1}\,\circ
\, t_{\,2}\,\circ \; \cdots \; \circ \,t_{\,k} \,)\; >\;0$. Then,
 \\$\mu\,(\,x_{\,n} \;-\; x \;,\; t_{\,1}\,\circ
\, t_{\,2}\,\circ \; \cdots \; \circ\, t_{\,k} \,) \;\,=\;\, \mu\,(\,
\sum\limits_{i\; = \;1}^k {\beta _{\,i}^{\,(\,n\,)} \,e_{\,i} }
\;-\; \sum\limits_{i\; = \;1}^k {\beta _{\,i} \,e_{\,i} } \;,\;\\{\hspace{9.5cm}} t_{\,1}\,\circ
\, t_{\,2}\,\circ \; \cdots \; \circ\, t_{\,k}
\,)\\  =\;\; \mu\,(\, \sum\limits_{i\; = \;1}^k
\,(\,{\beta _{\,i}^{\,(\,n\,)} \;-\; \beta_{\,i}\,) \,e_{\,i} }
\;,\; t_{\,1}\,\circ
\, t_{\,2}\,\circ \; \cdots \; \circ\, t_{\,k} \,) \\ \\ \ge \;\; \mu\,(\,(\,\beta
_{\,1}^{\,(\,n\,)} \;-\; \beta_{\,1} \,)\;e_{\,1} \;,\; t_{\,1} \,) \;\ast\; \cdots \;\ast\; \mu\,(\,(\,\beta _{\,k}^{\,(\,n\,)} \;-\;
\beta_{\,k} \,)\;e_{\,k} \;,\; t_{\,k} \,) \\\\  =
\;\;\mu\,(\,e_{\,1}
\;,\; \frac{t_{\,1}}{\,\psi (\,\beta _{\,1}^{\,(\,\,n\,)} \;-\; \beta_{\,1} \,)\,}
\;) \;\ast\; \cdots \;\ast\; \mu\,(\,e_{\,k}
\;,\; \frac{t_{\,k}}{\,\psi (\,\beta _{\,k}^{\,(\,\,n\,)} \;-\; \beta_{\,k} \,)\,}
\;)$ \\ Since \, $\mathop {\lim }\limits_{n\; \to \;\infty }
\;\frac{t_{\,i}}{\,\psi (\,\beta _{\,i}^{\,(\,\,n\,)} \;-\; \beta_{\,i} \,)\,}
\;\; = \;\infty $, we see that \\ \,$\mathop {\lim
}\limits_{n\; \to \;\infty }\;\mu\,(\,e_{\,i} \;,\;
\; \frac{t_{\,i}}{\,\psi (\,\beta _{\,i}^{\,(\,\,n\,)} \;-\; \beta_{\,i} \,)\,}
\;)
\;\,=\;\, 1$ \\ $\Longrightarrow\;\; \mathop {\lim }\limits_{n\; \to
\;\infty }\;\mu\,(\,x_{\,n} \;-\; x \;,\; t_{\,1}\,\circ
\, t_{\,2}\,\circ \; \cdots \; \circ \,t_{\,k}\,) \;\,\ge\;\, 1 \;\ast\;
\cdots \;\ast\; 1 \;=\; 1 $ \\
$\Longrightarrow\;\; \mathop {\lim }\limits_{n\; \to \;\infty
}\;\mu\,(\,x_{\,n} \;-\; x \;,\; t\,) \;\,=\;\, 1 $
  \, \, for all \,$t \;>\; 0$. \\\\ Again, we have
 \\$\nu\,(\,x_{\,n} \;-\; x\, \;,\; t_{\,1}\,\circ
\, t_{\,2}\,\circ \; \cdots \; \circ \,t_{\,k}\,) \;\,=\;\, \nu\,(\,
\sum\limits_{i\; = \;1}^k {\beta _{\,i}^{\,(\,n\,)} \,e_{\,i} }
\;-\; \sum\limits_{i\; = \;1}^k {\beta _{\,i} \,e_{\,i} } \;,\;\\{\hspace{9.5cm}} t_{\,1}\,\circ
\, t_{\,2}\,\circ \; \cdots \; \circ \,t_{\,k}
\,)\\  =\;\; \nu\,(\, \sum\limits_{i\; = \;1}^k
\,(\,{\beta _{\,i}^{\,(\,n\,)} \;-\; \beta_{\,i}\,) \,e_{\,i} }
\;,\; t_{\,1}\,\circ
\, t_{\,2}\,\circ \; \cdots \; \circ \,t_{\,k} \,) \\ \\ \ge \;\; \nu\,(\,(\,\beta
_{\,1}^{\,(\,n\,)} \;-\; \beta_{\,1} \,)\;e_{\,1} \;,\; t_{\,1} \,) \;\diamond\; \cdots \;\diamond\;  \nu\,(\,(\,\beta _{\,k}^{\,(\,n\,)} \;-\;
\beta_{\,k} \,)\;e_{\,k} \;,\; t_{\,k} \,) \\\\  =
\;\;\nu\,(\,e_{\,1}
\;,\; \frac{t_{\,1}}{\,\psi (\,\beta _{\,1}^{\,(\,\,n\,)} \;-\; \beta_{\,1} \,)\,}
\;) \;\diamond\; \cdots \;\diamond\; \nu\,(\,e_{\,k}
\;,\; \frac{t_{\,k}}{\,\psi (\,\beta _{\,k}^{\,(\,\,n\,)} \;-\; \beta_{\,k} \,)\,}
\;)$ \\ Since \, $\mathop {\lim }\limits_{n\; \to \;\infty }
\;\frac{t_{\,1}}{\,\psi (\,\beta _{\,1}^{\,(\,\,n\,)} \;-\; \beta_{\,1} \,)\,}
\;\; = \;\infty $, \, we see that \\ $\mathop {\lim
}\limits_{n\; \to \;\infty }\;\nu\,(\,e_{\,i} \;,\;
\; \frac{t_{\,1}}{\,\psi (\,\beta _{\,i}^{\,(\,\,n\,)} \;-\; \beta_{\,i} \,)\,}
\;)\;\,=\;\, 0 \;\; \forall \;i$
\\ $\Longrightarrow\;\; \mathop {\lim }\limits_{n\; \to
\;\infty }\;\nu\,(\,x_{\,n} \;-\; x \;,\; t_{\,1}\,\circ
\, t_{\,2}\,\circ \; \cdots \; \circ \,t_{\,k}\,) \;\,\le\;\, 0 \;\diamond\;
\cdots \;\diamond\; 0 \;=\; 0 $ \\
$\Longrightarrow\;\; \mathop {\lim }\limits_{n\; \to \;\infty
}\;\nu\,(\,x_{\,n} \;-\; x \;,\; t\,) \;\,=\;\, 0 \;\, \forall \; t \,>\,0$
 \\ \, Thus, we see that
\,$\{\,x_{\,n}\,\}_{\,n}$\, is an arbitrary Cauchy sequence that
converges to \,$x \;\,\varepsilon\;\, V$ \, hence the generalized IF$\psi$NLS \,$\left( {\,V\;,\;A^{\,\psi } \,} \right)$\, is complete .


\section{ Generalized Intuitionistic Fuzzy \\ $\psi$-Continuous Function}


\begin{Definition}
Let \,$\left( {\,U\;,\;A^{\,\psi } \,} \right)$\, and \,$(\,V \;,\;
B^{\,\psi }\,)$\, be two generalized IF$\psi$NLS
over the same field \,$F$. A mapping \,$f$\, from \,$\left( {\,U\;,\;A^{\,\psi } \,} \right)$\, to \,$(\,V \;,\;
B^{\,\psi }\,)$\, is said to be
\textbf{intuitionistic fuzzy continuous} \,$($\, or in short IFC
\,$)$ at \,$x_{\,0} \;\,\varepsilon\;\, U$, if for any given
\,$\varepsilon \;>\; 0$, \,$\alpha \;\,\varepsilon\;\, (\,0 \;,\;
1\,)$\, , \,$\exists \;\; \delta \;=\; \delta\,(\,\alpha \;,\;
\varepsilon\,) \;\,>\;\, 0 \;,\; \beta \;=\; \beta\,(\,\alpha \;,\;
\varepsilon\,) \;\,\varepsilon\;\, (\,0 \;,\; 1\,)$\, such that for
all \,$x \;\,\varepsilon\;\, U$, \\ ${\hspace{1.5cm}} \mu_{\,U}\,(\,x
\;-\; x_{\,0} \;,\; \delta\,) \;\,>\;\, \beta
\;\;\Longrightarrow\;\; \mu_{\,V}\,(\,f\,(\,x\,) \;-\;
f\,(\,x_{\,0}\,) \;,\; \varepsilon\,) \;\,>\;\, \alpha$
\hspace{0.5cm} and \\${\hspace{1.5cm}} \nu_{\,U}\,(\,x \;-\; x_{\,0}
\;,\; \delta\,) \;\,<\;\, 1 \;-\; \beta \;\;\Longrightarrow\;\;
\nu_{\,V}\,(\,f\,(\,x\,) \;-\; f\,(\,x_{\,0}\,) \;,\; \varepsilon\,)
\;\,<\;\, 1 \;-\; \alpha$. \\ If \,$f$\, is continuous at each point
of \,$U$ , $f$\, is said to be IFC on \,$U$ .
\end{Definition}

\begin{Definition}
A mapping \,$f$\, from \,$\left( {\,U\;,\;A^{\,\psi } \,} \right)$\, to \,$(\,V \;,\;
B^{\,\psi }\,)$\, is said to be \textbf{strongly intuitionistic fuzzy
continuous} \,$($\, or in short strongly IFC \,$)$ at \,$x_{\,0}
\;\,\varepsilon\;\, U$, if for any given \,$\varepsilon \;\,>\;\, 0$\, ,
\,$\exists \;\; \delta \;=\; \delta\,(\,\alpha \;,\; \varepsilon\,)
\;\,>\;\, 0
$\, such that for all \,$x \;\,\varepsilon\;\, U$, \\
${\hspace{2.5cm}} \mu\,_{\,V}\,(\,f\,(\,x\,) \;-\; f\,(\,x_{\,0}\,)
\;,\; \varepsilon\,) \;\,\ge\;\,\mu\,_{\,U}\,(\,x \;-\; x_{\,0} \;,\;
\delta\,)$ \hspace{0.5cm} and
\\${\hspace{2.5cm}} \nu\,_{\,V}\,(\,f\,(\,x\,) \;-\; f\,(\,x_{\,0}\,)
\;,\; \varepsilon\,) \;\,<\;\,\nu\,_{\,U}\,(\,x \;-\; x_{\,0} \;,\;
\delta\,)$ . \\ $f$\, is said to be strongly IFC on \,$U$\, if
\,$f$\, is strongly IFC at each point of \,$U$ .
\end{Definition}

\begin{Definition}

A mapping \,$f$\, from \,$\left( {\,U\;,\;A^{\,\psi } \,} \right)$\, to \,$(\,V \;,\;
B^{\,\psi }\,)$\, is said to be \textbf{sequentially intuitionistic fuzzy
continuous} \,$($\, or in short sequentially IFC \,$)$ at \,$x_{\,0}
\;\,\varepsilon\;\, U$, if for any sequence
\,$\{\,x_{\,n}\,\}_{\,n}$\, , \,$x_{\,n} \;\,\varepsilon\;\, U \;\;
\forall \;\; n$\, , \, with \,$x_{\,n} \;\longrightarrow\; x_{\,0}$
\, in \,$\left( {\,U\;,\;A^{\,\psi } \,} \right)$\, implies \,$f\,(\,x_{\,n}\,)
\;\longrightarrow\; f\,(\,x_{\,0}\,)$\, in \,$(\,V \;,\;
B^{\,\psi }\,)$ ,
that is ,for any given $r \,\in\,(\,0 \,,\,1\,)$, \,$t \;>\; 0$ \,
$\exists n_0\,\in\,N$ \\ ${\hspace{0.9cm}}
\mu_{\,U}\,(\,x_{\,n} \;-\; x_{\,0} \;,\; t\,) \;>\;
1 \;-\; r$\, and \,$\nu_{\,U}\,(\,x_{\,n} \;-\; x_{\,0} \;,\; t\,)
\;<\; r \;\;\forall \; n \; > \; n_{\,0}$ \\
$ \Longrightarrow\; \mu_{\,V}\,(\,f(\,x_{\,n}\,) \;-\; f(\,x_{\,0}\,) \;,\; t\,) \;>\;
1 \;-\; r$\, and \,$\nu_{\,V}\,(\,f(\,x_{\,n}\,) \;-\; f(\,x_{\,0}\,) \;,\; t\,)
\;<\;r \;\;\forall \; n \; > \; n_{\,0}$ \\ If \,$f$\, is sequentially IFC at each point of \,$U$\,
then \,$f$\, is said to be sequentially IFC on \,$U$ .
\end{Definition}

\begin{Theorem}
Let \,$f$\, be a mapping from $\left( {\,U\;,\;A^{\,\psi } \,} \right)$\, to
\,$(\,V \;,\; B^{\,\psi }\,)$. If \,$f$\, strongly IFC then it is
sequentially IFC.
\end{Theorem}

{\bf Proof.}$\;\;$
Let \,$f \;:\; (\,U \;,\;A^{\,\psi }\,) \;\longrightarrow\; (\,V \;,\;
B^{\,\psi }\,)$\, be strongly IFC on \,$U$\, and \,$x_{\,0}
\;\,\varepsilon\;\, U$. Then for each \,$\varepsilon \;>\; 0$\, ,
\,$\exists \;\; \delta \;=\; \delta\,(\,x_{\,0} \;,\; \varepsilon\,)
\;\,>\;\, 0$\, such that for all \,$x\;\,\varepsilon\;\,U$ , \\\\
${\hspace{2.5cm}} \mu\,_{\,V}\,(\,f\,(\,x\,) \;-\; f\,(\,x_{\,0}\,)
\;,\; \varepsilon\,) \;\,\ge\;\,\mu\,_{\,U}\,(\,x \;-\; x_{\,0} \;,\;
\delta\,)$ \hspace{0.5cm} and
\\${\hspace{2.5cm}} \nu\,_{\,V}\,(\,f\,(\,x\,) \;-\; f\,(\,x_{\,0}\,)
\;,\; \varepsilon\,) \;\,<\;\,\nu\,_{\,U}\,(\,x \;-\; x_{\,0} \;,\;
\delta\,)$ \\\\ Let \,$\{\,x_{\,n}\,\}_{\,n}$\, be a sequence in
\,$U$\, such that \,$x_{\,n} \;\longrightarrow\; x_{\,0}$  in the space $\left( {\,U\;,\;A^{\,\psi } \,} \right)$, that is
, for any given $r \,\in\,(\,0 \,,\,1\,)$, \,$t \;>\; 0$ \,
$\exists n_0\,\in\,N$ \, such that
$\;\mu\,_{\,U}\,(\,x_n\,-\,x\,,\,t\,)\;>\;1\,-\,r$\,\,and\,\, $\nu\,_{\,V}\,(\,x_n\,-\,x\,,\,t\,)\;<\;r \;\;\forall \; n \; > \; n_{\,0}$
  \\ Again, we see that \\ ${\hspace{2.5cm}}
\mu\,_{\,V}\,(\,f\,(\,x_{\,n}\,) \;-\; f\,(\,x_{\,0}\,) \;,\;
\varepsilon\,) \;\,\ge\;\,\mu\,_{\,U}\,(\,x_{\,n} \;-\; x_{\,0} \;,\;
\delta\,)$ \hspace{0.5cm} and
\\${\hspace{2.5cm}} \nu\,_{\,V}\,(\,f\,(\,x_{\,n}\,) \;-\; f\,(\,x_{\,0}\,)
\;,\; \varepsilon\,) \;\,<\;\,\nu\,_{\,U}\,(\,x_{\,n} \;-\; x_{\,0}
\;,\; \delta\,)$ \\ which implies that $\\ \mu\,_{\,V}\,(\,f(\,x_{\,n}\,) \;-\;
f(\,x_{\,0}\,) \;,\; \varepsilon\,) \;>\; 1 \;-\; r$\, and \,$\nu\,_{\,V}\,(\,f(\,x_{\,n}\,) \;-\;
f(\,x_{\,0}\,) \;,\; \varepsilon\,) \;<\;r \;\;\forall \; n \; > \; n_{\,0}$ \\ that is ,
\,$f\,(\,x_{\,n}\,) \;\longrightarrow\; f\,(\,x_{\,0}\,)$\, in
\,$(\,V \;,\; B^{\,\psi }\,)$. \, This completes the proof.

\begin{Theorem}
Let \,$f$\, be a mapping from the generalized IF$\psi$NLS
 \,$\left( {\,U\;,\;A^{\,\psi } \,} \right)$\, to
\,$(\,V \;,\; B^{\,\psi }\,)$. Then \,$f$\, is IFC on \,$U$\, if and only if
it is sequentially IFC on \,$U$ .
\end{Theorem}

{\bf Proof.}
The proof is same as the proof of the theorem 13 \cite{Samanta}

\begin{Theorem}
Let \,$f$\, be a mapping from the IFNLS \,$(\,U \;,\; A^{\,\psi }\,)$\, to
\,$(\,V \;,\; B^{\,\psi }\,)$\, and \,$D$\, be a compact subset of \,$U$ . If
\,$f$\, IFC on \,$U$\, then \,$f\,(\,D\,)$\, is a compact subset of
\,$V$ .

\end{Theorem}

{\bf Proof.}$\;\;$
Directly follows from the proof of the theorem 14 \cite{Samanta}

\end{document}